\documentclass[12pt]{amsart}
\usepackage{amsthm}
\usepackage{amsmath, mathrsfs}
\usepackage{amssymb}
\usepackage{hyperref}
\usepackage{enumerate}
\usepackage{latexsym,array}
\usepackage{amsfonts}
\usepackage{shadow}
\usepackage{tikz}
\usepackage{rotating}
\usepackage{tabularx}
\usepackage{graphicx}
\usepackage[section]{placeins}
\usepackage{mathtools}
\usepackage{url}

\newcommand{\bigH}{\mathcal{H}}
\newcommand{\bigA}{\mathcal{A}}
\newcommand{\bigC}{\mathcal{C}}
\newcommand{\bigD}{\mathcal{D}}

\newcommand{\bigK}{\mathcal{K}}
\newcommand{\bigR}{\mathcal{R}}

\newcommand{\bigB}{\mathcal{B}}

\newcommand{\bigO}{\mathcal{O}}
\newcommand{\bigT}{\mathcal{T}}
\newcommand{\bigV}{\mathcal{V}}
\newcommand{\bigJ}{\mathcal{J}}
\newcommand{\bigL}{\mathcal{L}}

\DeclareMathOperator{\Aut}{Aut}

\newtheorem{Pa}{Paper}[section]
\newtheorem{Tm}[Pa]{{\bf Theorem}}

\newtheorem{Dn}[Pa]{{\bf Definition}}
\newtheorem{Cy}[Pa]{{\bf Corollary}}

\newtheorem{Rk}[Pa]{{\bf Remark}}
\newtheorem{Pn}[Pa]{{\bf Proposition}}

\newenvironment{tightcenter}{%
  \setlength\topsep{0pt}
  \setlength\parskip{0pt}
  \begin{center}
}{%
  \end{center}
}

\tikzset{node distance=2cm, auto}

\date{}

\title[Dilation via non-commutative boundaries]
{Full Cuntz-Krieger dilations via non-commutative boundaries}


\author[A. Dor-On]{Adam Dor-On}
\address{Department of Pure Mathematics\\ University of Waterloo\\ Waterloo, ON\\ Canada}
\email{adoron@uwaterloo.ca}

\author[G. Salomon]{Guy Salomon}
\address{Department of Mathematics\\Technion --- Israel Institute of Technology\\Haifa \\Israel}
\email{guy.salomon@tx.technion.ac.il}

\thanks{The first author was partially supported by an Ontario Trillium scholarship. The second author was partially supported by the Clore Foundation.}

\begin{document}

\maketitle

\begin{tightcenter} \emph{\scriptsize 
We dedicate this article to the memory of Itzik Martziano, \\
who was an operator theorist, a brilliant colleague and a friend.}
\end{tightcenter}

\begin{abstract}
We apply Arveson's non-commutative boundary theory to dilate every Toeplitz-Cuntz-Krieger family of a directed graph $G$ to a full Cuntz-Krieger family for $G$. We do this by describing all representations of the Toeplitz algebra $\mathcal{T}(G)$ that have unique extension when restricted to the tensor algebra $\mathcal{T}_+(G)$. This yields an alternative proof to a result of Katsoulis and Kribs that the $C^*$-envelope of $\mathcal T_+(G)$ is the Cuntz-Krieger algebra $\mathcal O(G)$.

We then generalize our dilation results further, to the context of colored directed graphs, by investigating free products of operator algebras. These generalizations rely on results of independent interest on complete injectivity and a characterization of representations with the unique extension property for free products of operator algebras.

\end{abstract}

\section{Introduction}\label{sec:introduction}

Perhaps the simplest dilation result in operator theory is the dilation of an isometry to a unitary. If $V\in B(\bigH)$ is an isometry and $\Delta := I_{\bigH} - VV^*$, we may define a unitary $U$ on $K:= \bigH \oplus \bigH$ via
$$
U:= \begin{bmatrix}
V & \Delta \\
0 & V^*
\end{bmatrix}
$$
such that for any polynomial in a single variable $p\in \mathbb{C}[x]$ we have $p(V) = P_{\bigH} p(U)|_{\bigH}$ where $P_{\bigH}$ is the orthogonal projection onto the first summand of $\bigK = \bigH \oplus \bigH$. One of our goals in this paper is to generalize this dilation result to the free multivariable setting in the context of families of operators arising from directed graphs.

A {\em directed graph} $G$ is a quadruple $(V,E,s,r)$ consisting of a set $V$ of vertices, a set $E$ of edges and two maps $s,r : E\to V$, called the source map and range map, respectively. If $v=s(e)$ and $w=r(e)$ we say that $v$ emits $e$ and $w$ receives it. In this paper, we will assume that directed graphs are countable, meaning that both the sets $V$ and $E$ are countable. A directed graph is said to be {\em row-finite} if every vertex receives at most finitely many edges, and is {\em sourcelss} if every vertex receives at least one edge.

The tensor algebra $\bigT_+(G)$ and the $C^*$-algebras $\bigT(G)$ and $\bigO(G)$ associated to a directed graph $G$ have been studied by many authors. For instance \cite{Kat-Krib-Iso, Kribs-Power-Semigroupoid, Raeburn, Kumj-Pask-Raeburn-CK} and \cite{Cuntz-Krieger} to name but a few. We recommend \cite{Raeburn} and the references therein for additional details.

For a directed graph $G=(V,E,s,r)$ a {\em Toeplitz-Cuntz-Krieger $G$-family} $(P,S)$ is a set of mutually orthogonal projections $P:= \{P_v : v \in V\}$ and a set of partial isometries $S:= \{S_e : e \in E\}$ which satisfy the Toeplitz-Cuntz-Krieger relations:
\begin{enumerate}
\item[(I)] $S_e^*S_e=P_{s(e)}$ for every $e\in E$, and 
\item[(TCK)] $\sum_{e\in F} S_e S_e^* \leq P_{v}$ for every finite subset $F \subseteq r^{-1}(v)$.
\end{enumerate}
We say that $(P,S)$ is a {\em Cuntz-Krieger $G$-family} if, in addition, we have 
\begin{enumerate}
\item[(CK)]
$\sum_{r(e)=v} S_e S_e^* = P_v$ for every $v \in V$ with $0<|r^{-1}(v)|<\infty$.
\end{enumerate} 

The universal $C^*$-algebra $\bigT(G)$ generated by Toeplitz-Cuntz-Krieger $G$-families is called the {\em Toeplitz-Cuntz-Krieger algebra} of the graph $G$, and the universal $C^*$-algebra $\bigO(G)$ generated by Cuntz-Krieger families is called the {\em Cuntz-Krieger algebra} of the graph $G$. The {\em tensor algebra} $\bigT_+(G)$ is then just the norm-closed operator algebra generated by a universal Toeplitz-Cuntz-Krieger family and is a subalgebra of $\bigT(G)$. 

Due to their universal properties, $*$-representations of $\bigT(G)$ are in bijection with TCK $G$-families, and $*$-representations of $\bigO(G)$ are in bijection with CK $G$-families. Hence, we will often pass freely between these two points of view.

Toeplitz-Cuntz-Krieger families and Cuntz-Krieger families are easily seen to generalize the notions of an isometry and unitary respectively, by taking the graph with a single loop and a single vertex.

In our context of dilation of an isometry to a unitary, Popescu \cite[Proposition 2.6]{Popescu-isometric-dilations} proves that for a countable set $F$ and a row-isometry $V=(V_i)_{i\in F}$ on a space $\bigH$ there is a dilation to a row-unitary. In other words, this means that for any family of isometries $V_i : \bigH \rightarrow \bigH$ such that $\textsc{sot-}\sum_{i\in F}V_iV_i^* \leq I_{\bigH}$ there is a Hilbert space $\bigK$ containing $\bigH$, and isometries $U=(U_i)_{i\in F}$ on $\bigK$ such that $\textsc{sot-}\sum_{i\in F}U_iU_i^* = I_{\bigK}$, and for any polynomial $p\in \mathbb{C}\langle x_i \rangle_{i\in F}$ in non-commuting variables, we have $p(V) = P_{\bigH} p(U)|_{\bigH}$ where $P_{\bigH}$ is the projection from $\bigK$ to $\bigH$. In terms of graphs, this means that for a graph with a single vertex and $|F|$ loops, dilation of a TCK family to a CK family is possible, with the extra $\textsc{sot}$-convergence $\textsc{sot-}\sum_{i\in F}U_iU_i^* = I_{\bigK}$ when $F$ is infinite.

On the other hand, from \cite[Theorem 5.4]{Skalski-Zacharias-Wold} we see that if $G$ is row-finite and sourceless, then any TCK family has a CK dilation. More precisely for a row-finite sourceless graph $G=(V,E,s,r)$, if $(P,S)$ is a TCK family on $\bigH$, then there exists a CK family $(Q,T)$ on a larger space $\bigH$ such that for any polynomial $p\in \mathbb{C}\langle V,E \rangle$ in non-commuting variables we have $p(P,S) = P_{\bigH} p(Q,T)|_{\bigH}$. 

In order to put both of these results in the same context, we make the following definition. 

\begin{Dn}
Let $(P,S)$ be a Cuntz-Krieger family for a countable directed graph $G$. We say that $(P,S)$ is a \emph{full Cuntz-Krieger family} if 
\begin{enumerate}
\item[(CKF)]
$\textsc{sot-}\sum_{r(e)=v} S_e S_e^* = P_v$, for every $v\in V$ with $r^{-1}(v) \neq \emptyset$.
\end{enumerate}
\end{Dn}

We will show that in dilation theoretic terms, full CK families are the proper generalization of the notion of a unitary operator. More precisely, we will show that every TCK family has a full CK dilation, and that no non-trivial TCK dilations for full CK families are possible (see Corollary \ref{Cy:full-CK-dilation}).

We obtain our results by appealing to the non-commutative boundary theory introduced by Arveson \cite{Arveson_C-starI, Arveson_C-starII}. The notions of the $C^*$-envelope, also known as the non-commutative Shilov boundary, and the more delicate non-commutative Choquet boundary are very useful in operator algebras. A good instance of this appears in the recent work of Katsoulis and Ramsey, and Katsoulis on the Hao-Ng isomorphism problem \cite{Kats-Ram-semi, Kats-Hao-Ng}, where non-self-adjoint algebras and their $C^*$-envelopes play a prominent role. 

Let $\bigA$ be an operator algebra. We say that the pair $(\iota,\bigB)$ is a {\em $C^*$-cover} for a (not necessarily unital) operator algebra $\bigA$, if $\iota :\bigA \rightarrow \bigB$ is a completely isometric homomorphism and $C^*(\iota(\bigA)) = \bigB$. We will often identify $\bigA$ with its image $\iota(\bigA)$ under a given $C^*$-cover $(\iota,\bigB)$ for $\bigA$. We will call a linear map $\rho : \bigA \rightarrow B(\bigK)$ a \emph{representation} of $\bigA$ if it is a completely contractive homomorphism.

By \cite[Proposition 4.3.5]{BlecherLeMerdy} there is a unique, smallest $C^*$-cover for $\bigA$. This $C^*$-cover $(\iota, \bigB)$ is called the {\em $C^*$-envelope} of $\bigA$ and it satisfies the following universal property: given any other $C^*$-cover $(\iota', \bigB')$ for $\bigA$, there exists a (necessarily unique and surjective) $*$-homomorphism $\pi:\bigB' \to \bigB$, such that $\pi \circ \iota' = \iota$.

Characterizing the $C^*$-envelope of various operator structures was of use and intrigue to many authors, as can be seen for instance in \cite{Kak-Shalit, Dav-Ful-Kak, Kats-Ram-semi}. In \cite{Kats-Krib-Corresp}, Katsoulis and Kribs improve on the work in \cite{MuhlySolel98} and \cite[Theorem 5.3]{Fow-Muh-Rae} and show that the $C^*$-envelope of a tensor algebra associated to a general $C^*$-correspondence, is the Cuntz-Pimsner-Katsura algebra of the $C^*$-correspondence. In particular, the $C^*$-envelope of the tensor algebra $\bigT_+(G)$ is the Cuntz-Krieger algebra $\bigO(G)$ (this was also shown directly in \cite{KatsoulisKribs_Graphs}). We will provide an alternative proof for this fact on graph algebras in Theorem \ref{theorem:C-envelope-hyperrigidity}.

Suppose $\bigA$ is a \emph{unital} operator algebra generating a $C^*$-algebra $\bigB$. We say that a unital representation $\rho:\bigA \to B(\bigH)$ has the {\em unique extension property} if the only unital completely positive extension to $\bigB$ is a $*$-representation.

When $\bigA$ generates a $C^*$-algebra $\bigB$, Arveson defined boundary representations to be those irreducible $*$-representation of $\bigB$ whose restriction to $\bigA$ has the unique extension property. The collection of all these representations generalizes the classical notion of Choquet boundary for uniform algebras, and is therefore sometimes called the \emph{non-commutative Choquet boundary} for $\bigA$.

We say that $\bigA$ has the unique extension property in $\bigB$ if for any unital faithful $*$-representation $\pi : \bigB \rightarrow B(\bigH)$ we have that $\pi|_{\bigA}$ has the unique extension property. The unique extension property of $\bigA$ inside $C^*_{e}(\bigA)$ is equivalent to the notion of hyperrigidity, introduced by Arveson in \cite{Arv_HR}. Hyperrigidity and Arveson's hyperrigidity conjecture have been of interest to several authors recently. For instance, in \cite{Davidson-Kennedy-HR}, Davidson and Kennedy verify Arveson's hyperrigidity conjecture for commutative $C^*$-envelopes, and in the context of the Arveson-Douglas conjecture, Kennedy and Shalit show in \cite{Kennedy-Shalit-HR} that the essential normality of a $d$-tuple of operators satisfying homogeneous polynomial constraints is equivalent to the hyperrigidity of the $d$-tuple.

One of our main results is the classification of $*$-rep\-re\-sent\-a\-tion of $\bigT(G)$ that have the unique extension property when restricted to $\bigT_+(G)$. They turn out to coincide with those $*$-rep\-re\-sent\-a\-tions that are associated with full Cuntz-Krieger families (see Theorem \ref{Tm:non-annihi-no-uep}). This allows us to improve upon several known results and show that any TCK family dilates to a full CK family in the sense described above. Further applications of this result allows us to give a bijective correspondence between irreducible $*$-representations of $\bigT(G)$ that are not boundary, and ``gap" TCK families of the graph $G$ (see Corollary \ref{Cy:non-bdry-irred}), and a characterization of the unique extension property of $\bigT_+(G)$ in terms of the graph $G$ (see Theorem \ref{theorem:C-envelope-hyperrigidity}).

Trying to leverage our results to free products, we discuss some of the general theory of free products of operator algebras, and prove a joint unital completely positive extension theorem for free products of operator algebras amalgamated over any common $C^*$-algebra (see Theorem \ref{Tm:amalgamated_ucc}). Complete injectivity of amalgamated free products of $C^*$-algebras was shown by Armstrong, Dykema, Exel and Li \cite{ADER-free-prod}, and we are able to use our dilation techniques to generalize this to free products of operator algebras amalgamated over any common $C^*$-subalgebra (see Proposition \ref{Pn:CompleteInjectivity}).

In \cite[Theorem 5.3.21]{Dav-Ful-Kak} a gap in the proof of \cite[Theorem 3.1]{Duncan_FreeProd} was corrected, and it was shown that the amalgamated free product of $C^*$-envelopes is a $C^*$-cover for the amalgamated free product of operator algebras $\{\bigA_i\}_{i\in I}$. By \cite[Theorem 5.3.21]{Dav-Ful-Kak}, this $C^*$-cover turns out to be the $C^*$-envelope when each $\bigA_i$ has the unique extension property inside its $C^*$-envelope. We provide a characterization of representations with the unique extension property on amalgamated free products in Theorem \ref{tm:UEP-preservation-free-prod}, which yields these aforementioned results as well.

We apply our general results on free products to free products of operator algebras associated to graphs. These algebras have been investigated by Ara and Goodearl in \cite{Ara-Goodearl-separated} as $C^*$-algebras associated to separated graphs, and by Duncan \cite{Duncan_EdgeColored} as operator algebras associated to edge-colored directed graphs. We combine our results to prove that a full-CK dilation exists for any TCK family of a colored directed graph (see Corollary \ref{Cy:colored-dilation}), and to show that the free product of Cuntz-Krieger algebras is a $C^*$-cover for the free product of tensor graph algebras, which is the $C^*$-envelope when all graphs involved are row-finite (see Theorem \ref{Tm:freeHR}).

This paper has five sections including this introduction. In Section \ref{sec:ncb} we discuss some preliminary material on non-commutative boundary theory, especially in the non-unital context. In Section \ref{sec:Choquet-boundary} we describe a Wold decomposition for Toeplitz-Cuntz-Krieger families, and characterize representations whose restriction to the tensor algebra has the unique extension property. We use this to obtain our main dilation result and compute the $C^*$-envelope of graph tensor algebras. In Section \ref{sec:free-prod} we prove a joint extension theorem for free products of operator algebras amalgamated over a common $C^*$-algebra, along with the characterization of representations with the UEP. Finally, in Section \ref{sec:free-prod-graph}, we apply the results of Sections \ref{sec:Choquet-boundary} and \ref{sec:free-prod} to obtain a free product generalization, providing a free/colored version of our dilation result, and get our $C^*$-cover results for free products of graph tensor algebras.

\section{Non-commutative boundaries} \label{sec:ncb}

\subsection{$C^*$-envelopes, boundary representations and the unique extension property}\label{subsec:maximal}

Operator algebras can be given an axiomatic definition, as shown in \cite{BlecherRuanSinclair}. This means that they have an intrinsic operator structure that is preserved by any completely isometric homomorphism. We will survey the theory of non-commutative boundaries for operator algebras, and we refer the reader to \cite{Arveson_C-starI, Arveson_C-starII, Arveson-note, BlecherLeMerdy} for a more in-depth treatment of the theory.

For an operator algebra $\bigA$ generating a $C^*$-algebra $\bigB$, an ideal $\bigJ$ of $\bigB$ is called {\em a boundary ideal} for $\bigA$ if the quotient map $\bigB \to \bigB/\bigJ$ is completely isometric on $\bigA$. The largest boundary ideal $\bigJ_S(\bigA)$ of $\bigB$ is called {\em the Shilov ideal} of $\bigA$ in $\bigB$, and its importance in our context is that it gives a way to compute the $C^*$-envelope. Namely, the $C^*$-envelope of $\bigA$ is always isomorphic to $\bigB / \bigJ_S(\bigA)$.

When $\bigA$ is unital and $\pi : \bigB \to B(\bigH)$ is a unital $*$-representation such that $\pi|_{\bigA}$ has the unique extension property, every boundary ideal of $\bigA$ in $\bigB$ is contained in the kernel of $\pi$. The boundary theorem of Arveson in the separable case \cite{Arv08} and of Davidson and Kennedy \cite{Davidson-Kennedy} in general, then describes the Shilov ideal as the intersection of all kernels of boundary representations, providing yet another way to compute the $C^*$-envelope, via the non-commutative Choquet boundary.

For a (not-necessarily-unital) operator algebra $\bigA$ and a representation $\varphi: \bigA \rightarrow B(\bigH)$, a representation $\psi : \bigA \rightarrow B(\bigK)$ is said to {\em dilate} $\varphi$ if there is an isometry $V: \bigH \rightarrow \bigK$ such that for all $a\in \bigA$ we have $\varphi(a) = V^*\psi(a)V$. Since $V$ is an isometry, we can identify $\bigH \cong V(\bigH)$ as a subspace of $\bigK$, so that $\psi$ dilates $\varphi$ if and only if there is a larger Hilbert space $\bigK$ containing $\bigH$ such that for all $a\in \bigA$ we have that $\varphi(a) = P_{\bigH} \psi(a)|_{\bigH}$ where $P_{\bigH}$ is the projection onto $\bigH$.

In the case where $\bigA$ is unital, we say that a unital representation $\rho : \bigA \rightarrow B(\bigK)$ is {\em maximal} if whenever $\pi$ is a unital representation dilating $\rho$, then in fact $\pi = \rho \oplus \psi$ for some unital representation $\psi$. 

Building on ideas of Muhly and Solel \cite{MuhlySolel_boundary}, Dritschel and McCullough \cite[Theorem 1.1]{DritschelMcCullough} showed that a unital representation $\rho : \bigA \rightarrow B(\bigK)$ is maximal with respect to $\bigA$ if and only if it has the unique extension property with respect to $\bigA$. Dritschel and McCullough \cite[Theorem 1.2]{DritschelMcCullough} (see also \cite{Arveson-note}) then used this to show that every unital representation $\rho$ on $\bigA$ can be dilated to a \emph{maximal} unital representation $\pi$ on $\bigA$. This provided the first dilation-theoretic proof for the existence of the $C^*$-envelope.

\subsection{Non-commutative boundaries for non-unital algebras} \label{subsec:non-unital-tech}

We explain how to define the notions of maximality and the unique extension property for representations of \emph{not-necessarily-unital} operator algebras, in a way that yields the same theory as in the unital case.

If $\bigA\subseteq B(\bigH)$ is a \emph{non-unital} operator algebra generating a $C^*$-algebra $\bigB$, a theorem of Meyer \cite[Section 3]{Meyer} (see also \cite[Corollary 2.1.15]{BlecherLeMerdy}) states that every representation $\varphi: \bigA \rightarrow B(\bigK)$ extends to a unital representation $\varphi^1$ on the \emph{unitization} $\bigA^1=\bigA \oplus \mathbb{C}I_{\bigH}$ of $\bigA$  by specifying $\varphi^1(a+\lambda I_{\bigH}) = \varphi(a) + \lambda I_{\bigK}$. Meyer's theorem shows that $\bigA$ has a \emph{unique} (one-point) unitization, in the sense that if $(\iota, \bigB)$ is a $C^*$-cover for $\bigA$, and $\bigB \subseteq B(\bigH)$ is some faithful representation of $\bigB$, then the operator-algebraic structure on $\bigA^1 \cong \iota(\bigA) + \mathbb{C} 1_{\bigH}$ is independent of the $C^*$-cover and the faithful representation of $\bigB$.

Another important consequence is the not-necessarily-unital version of Arveson's extension theorem. More precisely, every representation $\varphi$ has a completely contractive and completely positive extension to $\bigB$ by first unitizing, applying the unital Arveson's extension theorem, and restricting back to the algebra. We record this as the following corollary.

\begin{Cy} \label{Cor:non-unital-arv}
Let $\bigA \subseteq B(\bigH)$ be an operator algebra generating a $C^*$-algebra $\bigB$, and let $\varphi : \bigA \rightarrow B(\bigK)$ be a representation of $\bigA$. Then there is a completely contractive and completely positive map $\widetilde{\varphi} : \bigB \rightarrow B(\bigK)$ such that $\hat{\varphi}|_{\bigA} = \varphi$.
\end{Cy}

Next, we discuss how to extend the notions of maximality and the unique extension property to not-necessarily-unital operator algebras.

\begin{Dn} \label{def:non-unital-UEP}
Let $\bigA \subseteq B(\bigH)$ be an operator algebra generating a $C^*$-algebra $\bigB$. Let $\rho : \bigA \rightarrow B(\bigK)$ be a representation.
\begin{enumerate}
\item
We say that $\rho$ has the \emph{unique extension property} (UEP for short) if every completely contractive and completely positive map $\pi: \bigB \rightarrow B(\bigK)$ extending $\rho$ is a $*$-representation.
\item
We say that $\rho$ is \emph{maximal} if whenever $\pi$ is a representation dilating $\rho$, then $\pi = \rho \oplus \psi$ for some representation $\psi$.
\end{enumerate}
\end{Dn}

It is important to note right away that in the case where $\bigA$ is unital and $\rho$ is unital, our definitions of maximality and the unique extension property reduce to the usual ones for unital representations of unital algebras. Indeed, for maximality, when $\pi:\bigA \to B(\bigL)$ is a representation dilating $\rho$, then $\pi(1)\pi(\cdot)|_{\pi(1)\bigL}$ is a unital dilation of $\rho$. So if $\rho$ is maximal only in the unital sense, then it must be a direct summand of $\pi(1)\pi(\cdot)|_{\pi(1)\bigL}$  and hence of $\pi$ itself, so that $\rho$ is maximal in the sense of Definition \ref{def:non-unital-UEP}. A similar argument using the unit works to show that $\rho$ has the unique extension property in the unital sense if and only if it has it in the sense of Definition \ref{def:non-unital-UEP}.

\begin{Rk} \em
When the maps in Definition \ref{def:non-unital-UEP} are not assumed to be multiplicative, there are instances where the UEP is satisfied vacuously. We thank Rapha\"el Clou\^atre for bringing these issues to our attention.

Indeed, Suppose $\bigA$ is a non-unital operator algebra containing a self-adjoint positive element $P$ and let $\rho : \bigA \rightarrow \bigB$ be a representation. The map $-\rho$ is completely contractive, but cannot be extended to a completely contractive completely positive map on $\bigB = C^*(\bigA)$, as $-\rho$ must send $P$ to $-P$. Hence, $-\rho$ vacuously has the UEP. Furthermore, when $\rho$ is \emph{not} maximal, the map $-\rho$ is a completely contractive map that admits a non-trivial completely contractive dilation, coming from the one for $\rho$. Hence, $-\rho$ is also not maximal. Thus, we see that if we drop the multiplicativity assumptions on our definitions above, the UEP and maximality would not be equivalent.
\end{Rk}

For the benefit of the reader, we provide a proof for the equivalence of the UEP with maximality in the not-necessarily-unital case.

\begin{Pn} \label{P:uep-max-non-unital}
Let $\bigA \subseteq B(\bigH)$ be an operator algebra generating a $C^*$-algebra $\bigB$, and let $\rho : \bigA \rightarrow B(\bigK)$ be a representation. Then $\rho$ has the UEP if and only if $\rho$ is maximal
\end{Pn}

\begin{proof}
Suppose that $\rho$ is maximal. Let $\phi$ be a completely positive and completely contractive extension of $\rho$ to $\bigB$. As necessary, using \cite[Proposition 2.2.1]{BrownOzawa} we can turn $\phi$ into a unital completely positive map $\phi^1$ on the unitization $\bigB^1$. Then use Stinespring's dilation theorem to get a dilation of $\phi^1$ to a unital $*$-representation $\sigma : \bigB^1 \rightarrow B(\bigL)$. By maximality we see that $\sigma|_{\bigA} = \rho \oplus \psi$ for some representation $\psi$ of $\bigA$, so that $\bigK$ is reducing for $\sigma|_{\bigB}$. Since the compression of $\sigma|_{\bigB}$ to $\bigK$ is $\phi$, we see that $\phi$ is multiplicative. Hence $\phi$ is unique among all completely positive completely contractive extensions. 

Conversely, if $\rho$ has the unique extension property, let $\phi : \bigA \rightarrow B(\bigL)$ be a representation dilating $\rho$. By Corollary \ref{Cor:non-unital-arv}, $\phi$ extends to a completely positive and completely contractive map $\widetilde{\phi} : \bigB \rightarrow B(\bigL)$. 
Compressing $\widetilde{\phi}$ to $\bigK$ gives a completely positive and completely contractive extension of $\rho$, which implies that $b \mapsto P_{\bigK} \widetilde{\phi}(b)|_{\bigK}$ is multiplicative. A standard use of Schwarz inequality \cite[Proposition 1.5.7]{BrownOzawa} (See the proof of \cite[Proposition 2.2]{Arveson-note}) then shows that $\bigK$ is reducing for $\widetilde{\phi}$. Hence, $\phi = \rho \oplus \psi'$ for some representation $\psi'$, and we see that $\rho$ is maximal.
\end{proof}

Consequentially, since maximality is an intrinsic property of the operator algebra, the unique extension property for representations does not depend on the choice of $C^*$-cover, even for non-unital operator algebras. We will often refer to this fact as the ``\emph{invariance of the UEP}". Next we show that when the operator algebra is non-unital, maximality (and hence the UEP) can be reduced to the unital case. 

When $\bigA$ is non-unital trouble can arise as the $0$ representation can fail to be maximal. Still, we may write $\rho = \rho_{nd} \oplus 0^{(\alpha)}$ where $\rho_{nd}$ is the non-degenerate part of $\rho$ and $0$ is the zero representation with some multiplicity $\alpha \geq 1$. In this case $\rho$ is maximal if and only if both $\rho_{nd}$ and $0$ are maximal. When $\bigA$ is unital, any issue with $0$ disappears since $0$ has UEP and is maximal automatically. Hence, we see that when $\bigA$ is unital, $\rho$ is maximal if and only if $\rho_{nd}$ is.

\begin{Pn} \label{P:unital-red}
Let $\bigA \subseteq B(\bigH)$ be a \emph{non-unital} operator algebra, and $\rho : \bigA \rightarrow B(\bigK)$ be a representation of $\bigA$. Then $\rho$ is maximal if and only if $\rho^1$ is maximal.
\end{Pn}

\begin{proof}
Suppose $\rho$ is maximal, and let $\psi : \bigA^1 \rightarrow B(\bigL)$ be a unital representation dilating $\rho^1$. Then clearly $\psi|_{\bigA}$ is a representation dilating $\rho$, so that $\bigK$ is reducing for $\psi(\bigA)$. Then clearly $\bigK$ is also reducing for $\psi(\bigA^1) = \psi(\bigA) + \mathbb{C}1_{\bigL} $, so that $\psi = \rho^1 \oplus \psi'$ for some unital representation $\psi'$. Hence $\rho^1$ is maximal.

Conversely, suppose that $\rho^1$ is maximal, and let $\psi$ be a representation dilating $\rho$. By Meyer's theorem we have a unital extension $\psi^1$ to a representation. Thus, we see that $\psi^1$ is a representation that dilates $\rho^1$, and by maximality of $\rho^1$, we have $\psi^1 = \rho^1 \oplus \psi'$ for some unital representation $\psi'$. If we restrict back to $\bigA$ we see that $\psi = \rho \oplus \psi'|_{\bigA}$, so that $\rho$ is maximal.
\end{proof}

Hence, from Proposition \ref{P:uep-max-non-unital} we see that when $\bigA$ is non-unital, a representation $\rho$ on $\bigA$ has the UEP if and only if its unitization $\rho^1$ has the UEP.

The $C^*$-envelope of a non-unital operator algebra can also be computed from the $C^*$-envelope of its unitization. More precisely, as the pair $(\kappa, C^*_e(\bigA))$ where $C^*_e(\bigA)$ is the $C^*$-subalgebra generated by $\kappa(\bigA)$ inside the $C^*$-envelope $(\kappa^1, C^*_e(\bigA^1))$ of the (unique) unitization $\bigA^1$ of $\bigA$. By the proof of \cite[Proposition 4.3.5]{BlecherLeMerdy} this $C^*$-envelope of an operator algebra $\bigA$ has the desired universal property, that for any $C^*$-cover $(\iota', \bigB')$ of $\bigA$, there exists a (necessarily unique and surjective) $*$-homomorphism $\pi:\bigB' \to C^*_e(\bigA)$, such that $\pi \circ \iota' = \kappa$.

As to representations with the UEP, when $\bigA$ is an operator algebra generating a $C^*$-algebra $\bigB$, through unitization the theorem of Dritschel and McCullough in the unital case shows that $C^*_e(\bigA)$ is again the image of a $*$-representation $\rho : \bigB \rightarrow B(\bigK)$ such that $\rho|_{\bigA}$ is completely isometric and has the unique extension property.

Let $\bigA$ be an operator algebra generating a $C^*$-algebra $\bigB$. We say that $\bigA$ has the unique extension property in $\bigB$ if for any faithful $*$-representation $\pi : \bigB \rightarrow B(\bigH)$ we have that $\pi|_{\bigA}$ has the unique extension property. By taking a direct sum of $\pi$ with a given $*$-representation of $\bigB$, it is easy to show that the faithfulness assumption can be dropped, and in particular, we must have that $\bigB \cong C^*_{e}(\bigA)$.

We will need the following easy consequence on the existence of a largest sub-representation with the UEP. Let $\phi: \bigB \rightarrow B(\bigH)$ be a completely contractive completely positive map on a $C^*$-algebra $\bigB$, and let $\bigK \subseteq \bigH$ be a reducing subspace for $\phi(\bigA)$. Let $\phi_{\bigK} : \bigB \rightarrow B(\bigK)$ denote the restriction $\phi_{\bigK}(b) = \phi(b)|_{\bigK}$.

\begin{Pn} \label{prop:largest-UEP}
Let $\bigA$ be an operator algebra generating a $C^*$-algebra $\bigB$ and let $\pi : \bigB \rightarrow B(\bigH)$ be a $*$-representation. Then there is a unique (perhaps trivial) largest reducing subspace $\bigK$ for $\pi$ such that $\pi_{\bigK}|_{\bigA}$ has the unique extension property.
\end{Pn}

\begin{proof}
If there is no such non-trivial reducing subspace, we take $\bigK = \{0\}$. Otherwise,
let $\mathcal C$ be the (non-empty) collection of non-trivial reducing subspaces $\bigL$ for $\pi$ such that $\pi_{\bigL} : \bigB \rightarrow B(\bigL)$ has the UEP when restricted to $\bigA$. Set $\bigK:=\bigvee _{\bigL \in \bigC} \bigL$. Since every $\bigL \in \bigC$ is reducing for $\pi$, we must have that $\bigK$ is reducing for $\pi$ as well.
It remains to show that $\pi_{\bigK}|_{\bigA}$ is maximal. So suppose that $\psi$ is a dilation of $\pi_{\bigK}|_{\bigA}$. Then by maximality of each $\pi_{\bigL}|_{\bigA}$ for $\bigL \in \bigC$, we see that $\bigL$ is reducing for $\psi$. Hence, so must $\bigK$ be reducing for $\psi$, and thus $\pi_{\bigK}$ must be maximal. 
\end{proof}

\section{Boundaries arising from directed graphs} \label{sec:Choquet-boundary}

Let $G = (V,E,s,r)$ be a countable directed graph. We will abuse terminology and call associated $*$-representations of either $\bigT(G)$ or $\bigO(G)$ ``Cuntz-Krieger" or ``full Cuntz-Krieger" if their associated TCK families are such. A (universal) TCK or CK family generating $\bigT(G)$ or $\bigO(G)$ (respectively) will usually be denoted by lowercase letters $(p,s)$.

There is a canonical $*$-representation of the Toeplitz-Cuntz-Krieger graph $C^*$-algebra which we now describe.
First, recall that a path in $G$ is a sequence of edges $\lambda =\mu_n \cdots \mu_1$ such  that $r(\mu_i) = s(\mu_{i+1})$, where we extend the range and source maps to apply for paths by specifying $r(\lambda) := r(\mu_n)$ and $s(\lambda) := s(\mu_1)$, and set $|\lambda|:=n$ for the length of the path; vertices are considered as paths of length $0$.  We use $E^{\bullet}$ to denote the collection of all paths in $G$ of finite length. 

Let $\bigH_G:= \ell^2(E^{\bullet})$ be the Hilbert space with canonical standard orthonormal basis $\{\xi_{\lambda}\}_{\lambda \in E^{\bullet}}$, we define a Toeplitz-Cuntz-Krieger family $(P,S)$ on $\bigH_G$ by specifying each operator on an orthonormal basis, that is, for each $v\in V$, $\mu \in E$ and $\lambda \in E^{\bullet}$ we define
\[
P_v(\xi_{\lambda}) = \begin{cases} 
\xi_{\lambda} & \text{if } r(\lambda) = v \\ 
0 & \text{if } r(\lambda) \neq v
\end{cases} \ \ \text{and} \ \
S_e(\xi_{\lambda}) = \begin{cases} 
\xi_{e \lambda} & \text{if } r(\lambda) = s(e) \\ 
0 & \text{if } r(\lambda) \neq s(e)
\end{cases}.
\]

For every $v\in V$, consider the subspace $\bigH_{G,v}:= \ell^2(s^{-1}(v))$ with its orthonormal basis
$\{ \xi_{\lambda} \}_{s(\lambda)=v}$. Clearly, $\bigH_{G,v}$ is reducing for $(P,S)$, so by the universal property of $\bigT(G)$ there exists a $*$-representation $\pi_v:\bigT(G) \to B(\bigH_{G,v})$ satisfying $\pi_v(p_w)=P_w|_{\bigH_{G,v}}$ for every $w \in V$ and  $\pi_v(s_e)=S_e|_{\bigH_{G,v}}$ for every $e \in E$. The next proposition is easily verified, and we omit its proof.

\begin{Pn}\label{Pn:H_G,v}
Let $\pi_v:\bigT(G) \to B(\bigH_{G,v})$ be the $*$-representation described above. Then the following hold:
\begin{enumerate}
\item[(a)] $\pi_v$ is irreducible,
\item[(b)] for every $w \neq v$ we have $\textsc{sot-}\sum_{r(e) = w} \pi_v (s_e s_e^*) = \pi_v(p_w)$, and
\item[(c)] $\pi_v(p_v)- \textsc{sot-}\sum_{r(e) = v}\pi_v (s_e s_e^*)$ is a rank $1$ projection.
\end{enumerate}
\end{Pn}

Toeplitz-Cuntz-Krieger families have the following useful version of the Wold decomposition. A slightly different Wold decomposition was given in \cite[Section 2]{Jury-Kribs} by Jury and Kribs under the assumption that the graphs have no sinks. Here we give a self-contained, and slightly more general version, that is tailored to our context. Let $V_r$ be the set of vertices $v \in V$ such that $r^{-1}(v) \neq \emptyset$. For a TCK family $(Q,T)$ and a reducing subspace $\bigK$ for it, we will denote $(Q,T)|_{\bigK} := (\{Q_v|_{\bigK}\},\{T_e|_{\bigK}\})$.

\begin{Tm}[Wold decomposition] \label{thm:wold-decomp}
Let $(Q,T)$ be a Toeplitz-Cuntz-Krieger family on a Hilbert space $\bigH$. For every $v \in V_r$, denote by $\alpha_v$ the dimension of the space $W_v:= (Q_v - \sum_{r(e)=v} T_e T_e^*)\bigH$. Then $(Q,T)$ is unitarily equivalent to
\[
\oplus_{v \in V_r} ((P,S)|_{\bigH_{G,v}})^{(\alpha_v)} \oplus (R,L)
\]
where $(R,L)$ is a full CK $G$-family. In addition, this representation is unique in the sense that if $(Q,T)$ is unitarily equivalent to 
\[
\oplus_{v \in V_r} ((P,S)|_{\bigH_{G,v}})^{(\alpha_v')} \oplus (R',L')
\]
where $(R',L')$ is a full CK $G$-family, then $\alpha_v'=\alpha_v$ for every $v \in V_r$, and $(R,L)$ is unitarily equivalent to $(R',L')$.
\end{Tm}

\begin{proof}
Uniqueness follows by Proposition \ref{Pn:H_G,v}. Indeed, as $(P,S)|_{\bigH_{G,v}}$ cannot be unitarily equivalent to $(P,S)|_{\bigH_{G,w}}$ for $w\neq v$ nor to restrictions to reducing subspaces for either full CK families $(R',L')$ or $(R,L)$. Thus, we must have that $((P,S)|_{\bigH_{G,v}})^{(\alpha_v)}$ is unitarily equivalent to $((P,S)|_{\bigH_{G,v}})^{(\alpha_v')}$ so that $\alpha_v = \alpha_v'$. Once this is established, restricting to the orthocomplement of the (reducing) subspaces associated with $\oplus_{v \in V_r} ((P,S)|_{\bigH_{G,v}})^{(\alpha_v)}$ and $\oplus_{v \in V_r} ((P,S)|_{\bigH_{G,v}})^{(\alpha_v')}$, we obtain a unitary equivalence between $(R,L)$ and $(R',L')$.

As for existence, fix $v\in V_r$, and denote $W_v =(Q_v - \sum_{r(e)=v}T_eT_e^*)\bigH$. Choose an orthonormal basis $\{\zeta_v^{(i)} \}$  for $W_v$, of cardinality $\alpha_v$, and for every $i$ set
\[
\bigH_{v,i}:={\rm{span}}\{T_{\lambda}\zeta_v^{(i)}: \lambda \in s^{-1}(v)\}.
\] 
We will show these subspace are reducing. Indeed, $\bigH_{v,i}$ is clearly invariant for the family $(Q,T)$. As for co-invariance, note that $T_{\mu}^* (T_{\lambda}\zeta_v^{(i)})$ is either $0$, a vector of the form $T_{\lambda'}\zeta_v^{(i)}$ for some path $\lambda'$, or a vector of the form $T_{\mu'}^*\zeta_v^{(i)}$ for some path $\mu'$ with $|\mu'|\geq 1$. As the two first cases immediately imply that $T_{\mu}^* (T_{\lambda}\zeta_v^{(i)}) \in \bigH_{v,i}$, we need to deal only with the third case.
To this end, write $\mu' = e_0 \mu''$ for some edge $e_0 \in E$ and a path $\mu''$ with $s(\mu'')=r(e_0)$. Note that if $e_0 \in s^{-1}(v)$, then 
$T_{e_0}^*(Q_v - \sum_{r(e)=v}T_eT_e^*)=T_{e_0}^* - T_{e_0}^*=0$, and otherwise,
$T_{e_0}^*(Q_v - \sum_{r(e)=v}T_eT_e^*)=0-0=0$.
Thus, in any case,
\[
T_{\mu'}^*\zeta_v^{(i)}    = T_{\mu'}^*(Q_v - \sum_{r(e)=v}T_eT_e^*)\zeta_v^{(i)}
                                        = T_{\mu''}^*T_{e_0}^*(Q_v - \sum_{r(e)=v}T_eT_e^*)\zeta_v^{(i)}
                                        =0.
\]

We next show simultaneously that for fixed $v\in V_r$ and $1 \leq i \leq \alpha_v$, the set $\{T_{\lambda}\zeta_v^{(i)}\}_{\lambda \in s^{-1}(v)}$ is an orthonormal family, and that the spaces $\bigH_{v,i}$ are pairwise orthogonal for all $v\in V_r$ and $1 \leq i \leq \alpha_v$. Our first step is to show that for two vertices $v,w\in V_r$, two indices $1\leq i \leq \alpha_v$ and $1 \leq j \leq \alpha_w$, and two paths $\lambda, \mu$ in $G$, if $\langle T_{\lambda} \zeta_v^{(i)} , T_{\mu} \zeta_w^{(j)} \rangle \neq 0$ then we must have $\lambda = \mu$.
Indeed, 
\[
\langle T_{\lambda} \zeta_v^{(i)} , T_{\mu} \zeta_w^{(j)} \rangle
=\Big\langle \big( (Q_w - \sum_{r(e)=w}T_eT_e^*) T_{\mu}^* T_{\lambda} (Q_v - \sum_{r(e)=v}T_eT_e^*) \big) \zeta_v^{(i)}, \zeta_w^{(j)} \Big \rangle.
\]
For $T_{\mu}^*T_{\lambda}$ to be non-zero, it must be either of the form $T_{\lambda'}$ where $\lambda = \mu \lambda'$, or $T_{\mu'}^*$ where $\mu = \lambda \mu'$. We deal with the first case, and the second is proven similarly. So assume $\lambda = \mu \lambda'$. If $|\lambda'| = 0$, then $\lambda = \mu$. If $|\lambda'| \geq 1$, write $\lambda' = e_0 \lambda''$. Then we have
\[
\big(Q_w - \sum_{e \in r^{-1}(w)}T_eT_e^*\big) T_{\mu}^* T_{\lambda} = T_{\lambda'} - T_{e_0}T_{e_0}^*T_{\lambda'} = 0
\]
which yields a contradiction. Thus, $\lambda = \mu$.

As a consequence of this, we see that $v = s(\lambda) = s(\mu) =w$. As $T_{\lambda}$ is an isometry on $P_v\bigH$, the assumption $\langle T_{\lambda} \zeta_v^{(i)} , T_{\lambda} \zeta_v^{(j)} \rangle \neq 0$ yields $i=j$ as well. We therefore must have that the sets $\{T_{\lambda}\zeta_v^{(i)}: \lambda\in s^{-1}(v)\}$ are orthonormal bases for the pairwise orthogonal reducing subspaces $\bigH_{v,i}$.

We next define unitaries $U_{v,i} : \bigH_{v,i} \rightarrow \bigH_{G,v}$ by mapping an orthonormal basis to an orthonormal basis $U_{v,i}: T_{\lambda}\zeta_v^{(i)} \mapsto \xi_{\lambda}$. Clearly $U_{v,i}$ intertwines $(Q,T)|_{\bigH_{v,i}}$ and $(P,S)|_{\bigH_{G,v}}$. Denote by $\bigK = (\oplus_{v\in V_r} \bigH_{v,i})^{\perp}$, we then have that $(Q,T)$ is unitarily equivalent to 
\[
\oplus_{v \in V_r} ((P,S)|_{\bigH_{G,v}})^{(\alpha_v)} \oplus (Q,T)|_{\bigK}.
\]
As a result $(R,L): = (Q,T)|_{\bigK}$ is a Topelitz-Cuntz-Krieger family such that for any $v\in V_r$ we have $R_v = \textsc{sot}-\sum_{r(e)=v}L_eL_e^*$. Hence, it is a full CK family.
\end{proof}

By rephrasing the previous proposition in terms of $*$-representations, we obtain the following corollary.

\begin{Cy} \label{cor:wold-decomp-rep}
Let $G$ be a directed graph, and let $\pi: \bigT(G) \rightarrow B(\bigH)$ be a $*$-representation. Then there are multiplicities $\{\alpha_v\}_{v\in V_r}$ such that $\pi$ is unitarily equivalent to the $*$-representation $\pi_s \oplus \pi_b$, where $\pi_s = \oplus_{v\in V_r}\pi_v^{(\alpha_v)}$ and $\pi_b$ is a full CK representation. In addition, this representation is unique in the sense that if $\pi$ is also unitarily equivalent to the $*$-representation $\pi_s' \oplus \pi_b'$, where $\pi_s' = \oplus_{v\in V_r}\pi_v^{(\alpha_v')}$ and $\pi_b'$ is a full CK representation, then $\alpha_v'= \alpha_v$ for every $v \in V_r$, and $\pi_b'$ is unitarily equivalent to $\pi_b$.
\end{Cy}

We next characterize those $*$-representations which have the unique extension property with respect to $\bigT_+(G)$.

\begin{Dn}
Let $G = (V,E,s,r)$ be a directed graph, and let $\pi : \bigT(G) \rightarrow B(\bigH)$ be a $*$-representation. We say that $v\in V_r$ is {\em singular} with respect to $\pi$ (or simply that $v$ is $\pi$-singular) if 
\[
\textsc{sot-}\sum_{r(e) = v} \pi (S_e S_e^*) \lneq \pi(P_v).
\]
\end{Dn}

Note that $\pi$ is a full CK representation of $\bigT(G)$, if and only if $\pi$ has no singular vertices.
 
\begin{Tm} \label{Tm:non-annihi-no-uep}
Suppose that $\pi : \bigT(G) \rightarrow B(\bigH)$ is a $*$-representation. The restriction $\pi|_{\bigT_+(G)}$ has the unique extension property if and only if $\pi$ is a full CK representation. 
\end{Tm}

\begin{proof}
Let $(p,s)$ be a generating TCK family for $\bigT(G)$. Suppose $G$ has a $\pi$-singular vertex $v$. If we assume towards contradiction that $\pi|_{\bigT_+(G)}$ has the unique extension property, then by \cite[Proposition 4.4]{Arv_HR}, so does the restriction of the infinite inflation $\pi^{(\infty)} : \bigT(G) \rightarrow B(\bigH^{(\infty)})$ to $\bigT_+(G)$. We therefore may assume without loss of generality that $\pi$ has infinite multiplicity. We will arrive at a contradiction by showing that $\pi|_{\bigT_+(G)}$ is not maximal. 

As $v$ is $\pi$-singular, and $\pi$ has infinite multiplicity, the projection $Q_v : = \pi(p_v) -\textsc{sot-} \sum_{r(e) = v} \pi (s_e s_e^*)$ is infinite dimensional. Thus, we may decompose $Q_v \bigH = \oplus_{r(e) = v}\bigH_e$ into infinite dimensional spaces $\bigH_e$ for each $e\in r^{-1}(v)$. We can then define for every $e\in r^{-1}(v)$ some isometry $W_e : P_{s(e)}\bigH_G \rightarrow \bigH_e$ where $\bigH_G$ is the Hilbert space $\ell^2(E^{\bullet})$ and $(P,S)$ is the associated TCK $G$-family. We moreover define a $*$-representation $\rho : \bigT(G) \rightarrow B(\bigH \oplus \bigH_G)$ by specifying a Toeplitz-Cuntz-Krieger $G$-family
\[
\rho(p_v)=
\begin{bmatrix}
\pi(p_v) & 0\\
0 & P_v
\end{bmatrix} \text{ for all $v\in V$}
\]
and
\[
\rho(s_e)=
\begin{cases}
~\begin{bmatrix}
\pi(s_e) & W_e \\
0 & 0
\end{bmatrix} &\text{if }r(e)=v, \quad \text{and}\\ 
&\\ 
~\begin{bmatrix}
\pi(s_e) & 0\\
0 & S_e
\end{bmatrix}
&\text{otherwise}.
\end{cases}
\]
We show this defines a Toeplitz-Cuntz-Krieger family. Clearly, we need to verify only those relations which involve edges in $r^{-1}(v)$. For every $e \in r^{-1}(v)$

\[
\begin{split}
\rho(s_e)^*\rho(s_e)&=
\begin{bmatrix}
\pi(s_e)^* & 0\\
W_e^* & 0
\end{bmatrix}
\cdot 
\begin{bmatrix}
\pi(s_e) & W_e \\
0 & 0
\end{bmatrix}\\
&=
\begin{bmatrix} 
\pi(p_{s(e)}) & \pi(s_e)^*W_e \\
W_e^*\pi(s_e) & P_{s(e)}
\end{bmatrix}.
\end{split}
\]
As the range of $W_e$ is orthogonal to that of $\pi(s_e)$, we see that $\pi(s_e)^*W_e = W_e^*\pi(s_e) = 0$, so 
\[
\rho(s_e)^*\rho(s_e) = \rho(p_{s(e)}),
\]
and condition (I) is verified. Next, for every finite subset $F \subseteq r^{-1}(v)$
\[
\begin{split}
\sum_{e \in F}\rho(s_e)\rho(s_e)^*
&= \sum_{e \in F}
\begin{bmatrix}
\pi(s_e) & W_e \\
0 & 0
\end{bmatrix}
\cdot 
\begin{bmatrix}
\pi(s_e)^* & 0\\
W_e^* & 0
\end{bmatrix}\\
&=
\sum_{e \in F} \begin{bmatrix}
\pi(s_es_e^*) + W_e W_e^* & 0\\
0 & 0
\end{bmatrix}\\
& \leq \begin{bmatrix}
\pi(p_v) & 0\\
0 & P_v
\end{bmatrix} = \rho(p_v).
\end{split}
\]
where the inequality is true since $\{ \pi(s_es_e^*) \} \cup \{ W_eW_e^* \}$ is a collection of pairwise orthogonal projections dominated by $\pi(p_v)$.
We therefore have shown condition (TCK), and we conclude that $\rho|_{\bigT_+(G)}$ is a well-defined representation which dilates $\pi|_{\bigT_+(G)}$ non-trivially. Hence $\pi|_{\bigT_+(G)}$ is not maximal.

For the converse, suppose that $G$ has no $\pi$-singular vertices. Let $\widetilde \rho: \bigT_+(G) \to B(\bigK)$ be a maximal dilation of $\pi|_{\bigT_+(G)}$, and let $\rho:\bigT(G) \to B(\bigK)$ be its extension to a $*$-representation. Denote
\[
\rho(p_v)=
\begin{bmatrix}
\pi(p_v) & X_v\\
Y_v & Z_v
\end{bmatrix}\quad \text{and} \quad
\rho(s_e)=
\begin{bmatrix}
\pi(s_e) & X_e\\
Y_e & Z_e
\end{bmatrix} 
\]
for all $v\in V$ and $e \in E$. We have that $X_v=0$ and $Y_v=0$ for all $v\in V$. Indeed, let $v\in V$, and $P: K \rightarrow H$ the orthogonal projection onto $H$, then
\[
\begin{split}
P \rho(p_v)^*(1-P) \rho(p_v) P &= P \rho(p_v) P - P \rho(p_v) P \rho(p_v) P\\
&=\pi(p_v) - \pi(p_v)\pi(p_v) = 0.
\end{split}
\]
and the $C^*$-identity implies $Y_v = (1-P) \rho(p_v) P = 0$. As $\rho(p_v)$ is self-adjoint, we have $X_v = 0$ as well. 

Next, for all $e\in E$ we have $p_{s(e)}=s_e^*s_e$, so
\[
\begin{bmatrix}
\pi(p_{s(e)}) & 0\\
0 & *
\end{bmatrix}
=\rho(s_e^*s_e)=\rho(s_e)^*\rho(s_e)=
\begin{bmatrix}
\pi(s_e)^*\pi(s_e)+Y_e^*Y_e & *\\
* & *
\end{bmatrix}
\]
which implies $Y_e=0$ for all $e\in E$.

Finally, let $e\in E$ and let $v=r(e)$. For every finite subset $F$ of $r^{-1}(v)$,
we have $p_{v} \geq \sum_{f \in F}s_fs_f^*$, so
\[
\begin{split}
\begin{bmatrix}
\pi(p_{v}) & 0\\
0 & *
\end{bmatrix}
&= \rho(p_{v})
\geq \sum_{f\in F} \rho(s_f)\rho(s_f)^*\\
&=\sum_{f\in F}
\begin{bmatrix}
\pi(s_f)\pi(s_f)^*+X_fX_f^* & *\\
* & *
\end{bmatrix}.
\end{split}
\]
In particular, by compressing this inequality to $\bigH$ we obtain
\[
\sum_{f\in F} \pi(s_f)\pi(s_f)^*+X_fX_f^* \leq \pi(p_{v})
\]
for every finite subset $F$ of $r^{-1}(v)$. Since $v$ is not $\pi$-singular, we must have that
\[
\sup_F \sum_{f\in F} \pi(s_f)\pi(s_f)^* =\textsc{sot-}\sum_{r(f)=v} \pi(s_f)\pi(s_f)^*=\pi(p_v).
\]
We therefore obtain that $X_f=0$ for all $f\in r^{-1}(v)$, and in particular $X_e=0$. Since $\bigT(G)$ is generated as a $C^*$-algebra by $\bigT_+(G)$, we must have that $\rho$ has $\pi$ as a direct summand, and hence $\rho|_{\bigT_+(G)}$ is a trivial dilation of $\pi|_{\bigT_+(G)}$.
\end{proof}
The previous theorem gives rise to two interesting corollaries. The first is a parametrization of those irreducible $*$-rep\-re\-sent\-a\-tions of $\bigT(G)$ which are not boundary representations with respect to $\bigT_+(G)$, and the second is the dilation of TCK families to full CK families.

\begin{Cy}\label{Cy:non-bdry-irred}
For every vertex $v\in V_r$, the $*$-representation $\pi_v : \bigT(G) \rightarrow B(\bigH_{G,v})$ is the unique irreducible $*$-representation (up to unitary equivalence) for which $v$ is $\pi$-singular, so that the irreducible $*$-representations of $\bigT(G)$ which are not boundary for $\bigT_+(G)$ are paramet\-rized by $V_r$.
\end{Cy}

\begin{proof}
If $\pi$ is an irreducible $*$-representation that lacks the unique extension property on $\bigT_+(G)$, then by Theorem \ref{Tm:non-annihi-no-uep} there exists $v\in V_r$ which is $\pi$-singular. By the Wold decomposition (Corollary \ref{cor:wold-decomp-rep}), up to a unitary equivalence, $\pi$ must have $\pi_v$ as a subrepresentation, and by irreducibility, $\pi$ is unitarily equivalent to $\pi_v$.
\end{proof}

\begin{Cy} \label{Cy:full-CK-dilation}
Let $G=(V,E,s,r)$ be a countable directed graph, and $(P,S)$ a TCK family on $\bigH$. Then there exists a full CK family $(Q,T)$ on a Hilbert space $\bigK$ containing $\bigH$, such that $f(P,S) = P_{\bigH} f(Q,T)|_{\bigH}$ for any polynomial $f \in \mathbb{C}\langle V,E \rangle$ in non-commuting variables.
\end{Cy}

\begin{proof}
Let $\pi_{P,S} : \bigT(G) \rightarrow B(\bigH)$ be the $*$-representation of $\bigT(G)$ associated to $(P,S)$. By \cite[Theorem 1.2]{DritschelMcCullough} we can dilate $\pi_{P,S}|_{\bigT_+(G)}$ to a maximal representation $\tau : \bigT_+(G) \rightarrow B(\bigK)$, and without loss of generality, $\bigH$ is a subspace of $\bigK$. Hence, $\tau$ is the restriction to $\bigT_+(G)$ of a $*$-representation $\rho: \bigT(G) \rightarrow B(\bigK)$ such that $\rho|_{\bigT_+(G)}$ has the unique extension property. Let $(Q,T)$ be the TCK family associated to $\rho$. By Theorem \ref{Tm:non-annihi-no-uep} $(Q,T)$ is a full CK family, and it dilates $(P,S)$ in the sense that for every polynomial $f\in \mathbb{C}\langle V, E \rangle$ we have that $f(P,S)= P_{\bigH}f(Q,T)|_{\bigH}$.
\end{proof}

Our next goal is to construct, for every directed graph $G$, faithful full CK representations of $\bigO(G)$.
We do this by constructing certain universal CK families arising from backward paths. This construction seems to originate from the work on atomic representations of free semigroup algebras from \cite{DP1999}, and was successfully used in \cite{KatsoulisKribs_Graphs} to show that the $C^*$-envelope of row-finite sourceless higher-rank graph tensor algebra is its higher rank graph $C^*$-algebra.

Let $E^{\infty} = \{ \ \lambda \ | \ \lambda = e_1 e_2 e_3 \cdots  \ , s(e_i) = r(e_{i+1}), \ e_i\in E \ \}$ be the collection of all backward infinite paths in $G$, and extend the range map to $E^{\infty}$ by setting $r(\lambda) = r(e_1)$ for $\lambda = e_1 e_2 e_3\cdots \in E^{\infty}$. Let $E^{<\infty}$ be the collection of all finite paths, including paths of length $0$, emanating from sources, and set $E^{\leq \infty} = E^{\infty} \cup E^{<\infty}$ as a disjoint union.  

For each vertex $v\in V$ fix an element $\mu_v \in E^{\leq \infty}$ with $r(\mu_v) = v$. For $i \in \mathbb{N}$ define $\mu_{v,i} = e_{v,1}...e_{v,i}$ as the $i$-th truncation of $\mu_v$, where if $|\mu_v| \leq i$, $\mu_{v,i} = \mu_v$. Denote $\bigH_{v,i}$ the Hilbert space with the orthonormal basis $\{ \ \xi_{\lambda \mu_{v,i}^{-1}} \ | \ \lambda \in E^{\bullet} \ \}$ where $\lambda \mu_{v,i}^{-1}$ corresponds to the equivalence class of reduced products determined by $\lambda$ and $\mu_{v,i}$. We set $\Gamma:= \{ \ \lambda \mu_{v,i}^{-1} \ | \ \lambda \in E^{\bullet}, \ i\in \mathbb{N}, \ v\in V \ \}$, and let $\bigH_b:= \ell^2(\Gamma)$ denote the Hilbert space with orthonormal basis $\{\xi_{\lambda \mu^{-1}}\}_{\lambda \mu^{-1} \in \Gamma}$, which is unitarily equivalent to $\oplus_{v\in V}\big[ \bigvee_{i\in \mathbb{N}} \bigH_{v,i} \big]$, where $\bigH_{v,i}$ is identified as a subspace of $\bigH_{v,i+1}$ since $\lambda \mu_{v,i}^{-1}$ is identified with $\lambda e_{v,i+1} \mu_{v,i+1}^{-1}$ whenever $|\mu_v| > i$ and with $\lambda \mu_{v,i}^{-1}$ when $|\mu_v| \leq i$. We define a TCK family $(Q,T)$ on $\bigH_b$ by specifying it on the orthonormal basis $\{\xi_{\lambda \mu^{-1}}\}_{\lambda \mu^{-1} \in \Gamma}$,
\[
Q_v(\xi_{\lambda \mu^{-1}}) = \begin{cases} 
\xi_{\lambda \mu^{-1}} & \text{if } r(\lambda) = v \\ 
0 & \text{if } r(\lambda) \neq v
\end{cases}, \ \
T_e(\xi_{\lambda \mu^{-1}}) = \begin{cases} 
\xi_{e \lambda \mu^{-1}} & \text{if } r(\lambda) = s(e) \\ 
0 & \text{if } r(\lambda) \neq s(e)
\end{cases}.
\]
It is easy to verify that $(Q,T)$ is a full CK family. Indeed, $\xi_{\lambda \mu_{v,i}^{-1}}$, with with $|\lambda| > 1$ is in the range of $T_e$ for $e\in E$ such that $\lambda = e \lambda'$, and each $\xi_{\mu_{v,i}^{-1}}$ where $s(\mu) = s(\mu_{v,i})$ is not a source, is in the range of $T_{e_{v,i+1}}$ as $\mu_{v,i}^{-1}$ is identified with $e_{v,i+1}\mu_{v,i+1}^{-1}$.

By construction, each $Q_v$ is non-zero for all $v\in V$, and we let $\rho_{\infty}$ be the $*$-representation of $\bigT(G)$ associated to $(Q,T)$ above. Moreover, by construction of $\rho_{\infty}$, for each $z\in \mathbb{T}$ we get a well-defined unitary $U_z : \bigH_b \rightarrow \bigH_b$ by specifying $U_z(\xi_{\lambda \mu^{-1}}) = z^{|\mu| - |\lambda|}\xi_{\lambda \mu^{-1}}$. It is then easy to see that $U_zQ_v U_z^* =Q_v$ and $U_zT_eU_z^* = zT_e$ so we may define a gauge action $\alpha : \mathbb{T} \rightarrow \Aut(C^*(Q,T))$ via $\alpha_z(A) = U_zAU_z^*$, so that $\alpha_z(Q_v) = Q_v$ and $\alpha_z(T_e) = zT_e$. Hence, by the gauge invariant uniqueness theorem $\rho_{\infty}$ is injective. The advantage of this construction is that it produces a space which has a natural action of $\mathbb{T}$ on it.

\begin{Rk} \em
One may form a full CK family on $\ell^2(E^{\leq \infty})$ in a similar way. However, by \cite[Theorem 1.2]{Szy-Gen} this representation will fail to be injective when the graph has a vertex-simple cycle with no entry.
\end{Rk}

Let $\bigJ(G)$ denote the kernel of the quotient $q: \bigT(G) \rightarrow \bigO(G)$. Evidently, $\bigJ(G)$ is the ideal of $\bigT(G)$ generated by terms of the form $p_v - \sum_{r(e)=v}s_es_e^*$ for vertices $v$ with $0 < |r^{-1}(v)|< \infty$. In \cite[Theorem 3.3]{Kakariadis-Tensor}, Kakariadis showed that $\bigT_+(G)$ has the unique extension property in $\bigO(G)$ when $G$ is row-finite. We provide the proof for this statement along with its converse, and the computation of the $C^*$-envelope.

\begin{Tm} \label{theorem:C-envelope-hyperrigidity}
Let $G = (V,E,s,r)$ be a directed graph, and let $q : \bigT(G) \rightarrow \bigO(G)$ be the natural quotient map. Then $q$ is completely isometric on $\bigT_+(G)$, and $C^*_e(\bigT_+(G)) \cong \bigO(G)$.

Moreover, $\bigT_+(G)$ has the unique extension property in $\bigO(G)$ if and only if $G$ is row-finite. 
\end{Tm}

\begin{proof}
Let $(p,s)$ be a TCK family such that its associated $*$-represen\-tation $\pi_{p,s}$ is faithful. By \cite[Theorem 1.1]{DritschelMcCullough} we know that $\pi_{p,s}|_{\bigT_+(G)}$ has a maximal dilation, so let $\tau : \bigT(G) \rightarrow B(\bigK)$ be a $*$-representation such that $\tau|_{\bigT_+(G)}$ is the dilation of $\pi_{p,s}|_{\bigT_+(G)}$, so that it is completely isometric and with the unique extension property. By Theorem \ref{Tm:non-annihi-no-uep}, $\tau$ is a full CK representation, and hence annihilates the Cuntz-Krieger ideal $\bigJ(G)$. Hence, $\tau$ must factor through the quotient map $q : \bigT(G) \rightarrow \bigO(G)$ by $\bigJ(G)$, and we have that $q$ is completely isometric on $\bigT_+(G)$.

Next, we show that if $G$ is row-finite then $\bigT_+(G)$ has the unique extension property in $\bigO(G)$ via $q$. By Theorem \ref{Tm:non-annihi-no-uep}, we see that every $*$-representation of $\bigT(G)$ that annihilates $\bigJ(G)$ has the unique extension property when restricted to $\bigT_+(G)$ inside $\bigT(G)$. Since $q$ is completely isometric on $\bigT_+(G)$, by invariance of the unique extension property, we see that every $*$-representation of $\bigO(G)$ has unique extension property when restricted to $\bigT_+(G)$ inside $\bigO(G)$. By \cite[Theorem 1.1]{DritschelMcCullough} we know that the $C^*$-envelope of $\bigT_+(G)$ is the image under the direct sum of all $*$-representations of $\bigO(G)$ with the unique extension property, so that $C^*_e(\bigT_+(G)) \cong \bigO(G)$ when $G$ is row-finite, as all $*$-representations of $\bigO(G)$ have the unique extension property when restricted to $\bigT_+(G)$ inside $\bigO(G)$.

Otherwise, if $G$ is not row-finite, we have that $\rho_{\infty} \circ q$ is a full CK representation with kernel $\bigJ(G)$, and hence again by invariance of the unique extension property and Theorem \ref{Tm:non-annihi-no-uep} we have that $\rho_{\infty}$ has the unique extension property on $\bigT_+(G)$. Hence, since $\rho_{\infty}$ is faithful, we still have that $C^*_e(\bigT_+(G)) \cong \bigO(G)$.

For the converse of the second part of the statement, suppose that $G$ is not row-finite, and let $v\in V$ be an infinite receiver. Then $\pi_v$ annihilates $\bigJ(G)$, so we may consider the induced $*$-representation $\dot{\pi}_v : \bigO(G) \rightarrow B(\bigH_{G,v})$. By Corollary \ref{Cy:non-bdry-irred}, we see that $\pi_v$ does not have the unique extension property when restricted to $\bigT_+(G)$, so that by invariance of the unique extension property, $\dot{\pi}_v$ does not have the unique extension property when restricted to $\bigT_+(G)$. Thus, $\bigT_+(G)$ does not have the unique extension property in $\bigO(G)$.
\end{proof}

\begin{Rk}\em
In Theorem \ref{Tm:non-annihi-no-uep}, and Corollaries \ref{Cy:non-bdry-irred} and \ref{Cy:full-CK-dilation} we avoided the use of a uniqueness theorem for $\bigO(G)$. This is also true for the computation of the $C^*$-envelope in Theorem \ref{theorem:C-envelope-hyperrigidity} when $G$ is row-finite. A uniqueness theorem was needed only for the computation of the $C^*$-envelope when the graph is non-row-finite.
\end{Rk}


\section{Free products and the unique extension property} \label{sec:free-prod}

Consider the category of unital $\mathbb{C}$-algebras (with unital homomorphisms as morphisms). Let $\{\bigA_i\}_{i \in I}$ be a family of unital $\mathbb{C}$-algebras and let $\bigD$ be a common unital subalgebra with $1_{\bigA_i} \in \bigD$, let $\iota_i:\bigD \to \bigA_i$ denote the natural embeddings. Pushouts in this category are known to exist, and are called {\em free product of $\{\bigA_i\}_{i \in I}$ amalgamated over the common subalgebra $\mathcal D$}, denoted by $\underset{\bigD}{*} \bigA_i$. We recall the details briefly.

We let $* \bigA_i$ (with no $\mathcal D$) denote the vector space spanned by formal expressions $a_1 * \cdots * a_n$ where $a_1 \in \bigA_{i_1}, \dots, a_n \in \bigA_{i_n}$ such that $i_1 \neq i_2 \neq \cdots \neq i_n$ and $n\geq 1$, where this expression behaves multilinearly, and we define multiplication of two such expressions, where $b_1 \in \bigA_{j_1},\dots, b_m \in \bigA_{j_m}$ with $j_1 \neq \cdots \neq j_m$ via
\[
(a_1 * \cdots * a_n) \cdot (b_1 * \cdots * b_m) = \begin{cases} 
a_1 *\cdots * a_n * b_1 * \cdots * b_m, & \text{if } i_n \neq j_1, \\ 
a_1 *\cdots * (a_n \cdot b_1) * \cdots* b_m, & \text{if } i_n = j_1.
\end{cases}
\]
With this multiplication, $* \bigA_i$ becomes a $\mathbb{C}$-algebra generated by $\{\bigA_i\}_{i\in I}$. Next, we identify the different copies of $\bigD$ by taking a quotient by the ideal $\langle \iota_i(d)- \iota_j(d) \rangle_{i,j \in I, d\in \bigD}$. This quotient, which is denoted by $\underset{\bigD}{*} \bigA_i$, has the following universal property. If $\mathcal B$ is another unital $\mathbb{C}$-algebras with unital $\mathbb{C}$-homomorphisms $\psi_i :\bigA_i \to \bigB$ which agree on $\bigD$, then there is a unital $\mathbb{C}$-homomorphism $\underset{\bigD}{*} \psi_i: \underset{\bigD}{*} \bigA_i \to B$ extending each $\psi_i$ on $\bigA_i$.

Next, we construct free products in the category of unital operator algebras with unital completely contractive homomorphisms. Let $\{\bigA_i\}_{i \in I}$ be a family of unital operator algebras with $\bigD$ a common unital operator subalgebra with $1_{\bigA_i} = 1_{\bigD}$ for all $i\in I$, and let $\underset{\bigD}{*}\bigA_i$ denote the free product of $\{\bigA_i\}_{i \in I}$ amalgamated over $\bigD$ in the larger category of $\mathbb C$-algebras. We define matricial $n$-semi-norms
\[
\|a\|_{n} := \sup \left\|\underset{\bigD}{*} \psi_i ^{(n)}(a) \right\|_{B(\bigH)},  \quad \forall n\in \mathbb N,~a \in M_{n}\left(\underset{\bigD}{*}\bigA_i\right)
\]
where the supremum is taken over all families $\{ \psi_{i}:\bigA_i \to B(\bigH) \}_{i \in I}$ of unital completely contractive homomorphisms that agree on $\mathcal D$ and over a Hilbert space $\bigH$ of large enough cardinality. It follows that $\bigJ=\{a\in\underset{\bigD}{*}\bigA_i:\|a\|=0\}$ is a two-sided ideal, and we denote the matricial $n$-norms on the quotient by $\|\cdot\|_{n}$ as well. The matricial $n$-norms $\| \cdot\|_{n}$ then define an operator-algebraic structure on the completion ${\underset{\bigD}{\hat{*}}}\bigA_i$ of $\underset{\bigD}{*}\bigA_i/\bigJ$, by the Blecher-Ruan-Sinclair theorem \cite{BlecherRuanSinclair}. Furthermore, by construction there are unital completely contractive homomorphisms $\iota_j : \bigA_j \rightarrow {\underset{\bigD}{\hat{*}}}\bigA_i$. We will later see that each $\iota_j$ is in fact completely isometric, so that each $\bigA_j$ can be thought of as an operator subalgebra of ${\underset{\bigD}{\hat{*}}}\bigA_i$ via $\iota_j$.

The operator algebra ${\underset{\bigD}{\hat{*}}}\bigA_i$ is called the free product of $\{\bigA_i\}_{i\in I}$ amalgamated over the common operator subalgebra $\bigD$, and has the following universal property by construction: for any unital operator algebra $\bigB$ and $\psi_i:\bigA_i \to \bigB$ unital completely contractive homomorphisms which agree on $\bigD$, there exists a completely contractive homomorphism $\psi:= {\underset{\bigD}{{*}}}\psi_i$ from ${\underset{\bigD}{\hat{*}}}\bigA_i$ into $\bigB$ such that $\psi_i = \psi \circ \iota$.

We now provide a joint completely contractive extension result for free products of operator algebras amalgamated over a common $C^*$-subalgebra. Our proof is an adaptation of a proof given by Ozawa in \cite[Theorem 15]{Connes-embedding-Ozawa} in the case of amalgamation over the complex numbers. Recall that whenever $\bigD$ is a subalgebra of an operator algebra $\bigA$, a completely contractive map $\phi:\bigA \to B(\bigH)$ is said to be a {\em $\bigD$-bimodule} map if $\phi(a_1da_2)=\phi(a_1)\phi(d)\phi(a_2)$ for every $a_1,a_2 \in \bigA$ and $d \in \bigD$ . By \cite[Proposition 1.5.7]{BrownOzawa} the restriction of a completely contractive map $\phi$ to a $C^*$-subalgebra $\bigD$ is multiplicative if and only if $\phi$ is a $\bigD$-bimodule map.

\begin{Tm}\label{Tm:amalgamated_ucc}
Let $\{\bigA_i\}_{i\in I}$ be a family of unital operator algebras containing a common unital $C^*$-algebra $\bigD$ with $1_{\bigA_i} = 1_{\bigD}$, and let $\phi_i : \bigA_i \rightarrow B(\bigH)$ be unital completely contractive $\bigD$-bimodule maps that agree on $\bigD$. Then there exists a unital completely contractive map $\phi : \underset{\bigD}{\hat{*}}\bigA_i \rightarrow B(\bigH)$ such that $\phi \circ \iota_i = \phi_i$ for all $i \in I$.
\end{Tm}

\begin{proof}
We construct multiplicative dilations of $\phi_i$ which agree on $\bigD$, so that the compression of their free product to $\bigH$ yields a unital completely contractive joint extension $\phi$ as in the statement of the theorem.

First, we set $\bigH_{1} := \bigH$ and $\phi_i^{(1)}:=\phi_i$. By the Arveson-Stinespring dilation theorem we may dilate each of these to a completely contractive homomorphism from $\bigA_i$ to $B(\bigH_{1} \oplus \bigK^{(i)}_1)$. Denote by $\rho_i^{(1)}$ this dilation of $\phi_i$, so that $\phi_i(a) = P_{\bigH_1}\rho_i^{(1)}(a)|_{\bigH_1}$. We note that since $\bigD$ is a $C^*$-algebra, and each $\phi_i^{(1)}$ is multiplicative on $\bigD$, the space $\bigK^{(i)}_1$ is reducing for $\rho^{(1)}_i|_{\bigD}$, so that for all $d\in \bigD$ we have
\[
\rho^{(1)}_i(d) = \phi_i^{(1)}(d) \oplus P_{\bigK^{(i)}_1}\rho_i^{(1)}(d)|_{\bigK^{(i)}_1}.
\]

Now suppose we have a sequence of subspaces
\[
\bigH_1 \subseteq \bigH_2 \subseteq \cdots \subseteq \bigH_n
\]
such that for all $1\leq m \leq n$ we have unital completely contractive $\bigD$-bimodule maps $\phi_i^{(m)}: \bigA_i \rightarrow B(\bigH_m)$ that agree on $\bigD$, along with representations $\rho_i^{(m)}: \bigA_i \rightarrow B(\bigH_m \oplus \bigK^{(i)}_m)$ that dilate each $\phi_i^{(m)}$, so that for all $d\in \bigD$ we have 
\[
\rho^{(m)}_i(d) = \phi_i^{(m)}(d) \oplus P_{\bigK^{(i)}_m}\rho_i^{(m)}(d)|_{\bigK^{(i)}_m}
\]
and that $\rho_i^{(m+1)}$ and $\phi^{(m+1)}_i$ have $\rho_i^{(m)}$ as a direct summand for every $i\in I$ and $1 \leq m < n$.

Denote by $\bigH_{n+1} = \bigH_{n} \oplus \bigoplus_{i\in I} \bigK^{(i)}_n$. Fix $i\in I$, and consider the map $\tau_i: \bigD \rightarrow B(\bigK^{(i)}_n)$ given by $\tau_i(d) = P_{\bigK^{(i)}_n}\rho_i^{(n)}(d) |_{\bigK^{(i)}_n}$. By applying Arveson's extension theorem, followed by a restriction, for any $j \in I$ distinct from $i$, we may extend $\tau_i$ to a unital completely contractive map $\sigma_{ji} : \bigA_j \rightarrow B(\bigK^{(i)}_n)$. We define for all $j\in I$,
\begin{equation} \label{dilation-form}
\phi_j^{(n+1)}:= \rho_j^{(n)} \oplus \bigoplus_{j \neq i \in I} \sigma_{ji} : \bigA_j \rightarrow B(\bigH_{n+1}),
\end{equation}
so that $\rho_j^{(n)}$ is a direct summand of $\phi_j^{(n+1)}$. We then have for every $d\in \bigD$ that
\[
\phi_j^{(n+1)}(d) = \phi_j^{(n)}(d) \oplus \bigoplus_{i\in I} P_{\bigK^{(i)}_{n}}\rho_i^{(n)}(d)|_{\bigK^{(i)}_{n}}.
\]
Hence, since the maps $\{\phi_i^{(n)}\}_{i\in I}$ all agree on $\bigD$, we have that the maps $\{\phi_i^{(n+1)}\}_{i\in I}$ all agree on $\bigD$.

We use Arveson's extension, Stinespring's theorem and the special form of $\phi_j^{(n+1)}$ in equation \eqref{dilation-form} to obtain a \emph{multiplicative} unital completely contractive map $\rho_j^{(n+1)}: \bigA_j \rightarrow B(\bigH_{n+1} \oplus \bigK^{(j)}_{n+1})$ dilating $\phi_j^{(n+1)}$ such that each $\rho_j^{(n)}$ is still a direct summand of $\rho_j^{(n+1)}$. Hence, we then get that $\rho_j^{(n+1)}(a)|_{\bigH_{m} \oplus \bigK^{(j)}_{m}} = \rho_j^{(m)}(a)$ for all $1 \leq m \leq n$, and that for all $d\in \bigD$,
\[
\rho^{(n+1)}_i(d) = \phi_i^{(n+1)}(d) \oplus P_{\bigK^{(j)}_{n+1}}\rho_i^{(n+1)}(d)|_{\bigK^{(j)}_{n+1}}.
\]

Since for each $j\in I$ and $n \in \mathbb{N}$ we have $\bigH_{n} \subseteq \bigH_{n} \oplus \bigK^{(j)}_{n} \subseteq \bigH_{n+1}$, we may define a multiplicative unital completely contractive map $\rho_j: \bigA_j \rightarrow B(\bigK)$ on the inductive limit of Hilbert spaces $\bigK = \bigvee_{n\in \mathbb{N}} \bigH_{n}$ by specifying $\rho_j(a)h = \rho^{(n)}_j(a)h$ for $h\in \bigH_{n} \oplus \bigK^{(j)}_{n}$. These maps then agree on $\bigD$, since for $h\in \bigH_{n}$ we have that
\[
\begin{split}
\rho_i(d)h &= \rho_i^{(n+1)}(d)h \\
               &= \big( \phi_i^{(n+1)}(d) \oplus P_{\bigK^{(i)}_{n+1}}\rho_i^{(n+1)}(d)|_{\bigK^{(i)}_{n+1}} \big ) h\\  
               &= \phi_i^{(n+1)}(d) h
\end{split}
\]
and as the maps $\{\phi_i^{(n+1)}\}_{i\in I}$ all agree on $\bigD$ and the union of $\bigH_{n}$ is dense in $\bigK$, we have that $\rho_j(d) = \rho_i(d)$ for all $i\neq j$ in $I$.

Hence, we may form the free product $\rho:= \underset{\bigD}{*}\rho_i : \underset{\bigD}{\hat{*}}\bigA_i \rightarrow B(\bigK)$ which satisfies $\rho \circ \iota_i = \rho_i$, and the compression of $\rho$ to $\bigH$ would yield a joint unital completely contractive extension $\phi$ as in the statement of the theorem.
\end{proof}

From the above it is easy to see that if $\{\bigA_i\}_{i\in I}$ is a family of unital operator algebras containing a common unital $C^*$-algebra $\bigD$ with $1_{\bigA_i} = 1_{\bigD}$, then for each $i\in I$ the map $\iota_i : \bigA_i \rightarrow {\underset{\bigD}{\hat{*}}}\bigA_i$ is completely isometric. Hence, we will identify $\bigA_i$ as a unital operator subalgebra of $\underset{\bigD}{\hat{*}}\bigA_i$ via $\iota_i$

For the proof, fix $j \in I$ and let $\phi_j : \bigA_j \rightarrow B(\bigH)$ be a unital completely isometric homomorphism. We may then restrict it $\bigD$ and use Arveson's extension theorem to extend to unital completely contractive maps $\phi_i : \bigA_i \rightarrow B(\bigH)$ for $i\in I$ such that $i\neq j$. By Theorem \ref{Tm:amalgamated_ucc} there is a joint unital completely contractive map $\phi : {\underset{\bigD}{\hat{*}}}\bigA_i \rightarrow B(\bigH)$ which we may then dilate to a multiplicative map $\rho : {\underset{\bigD}{\hat{*}}}\bigA_i \rightarrow B(\bigK)$. However, the compression of $\phi \circ \iota_j$ to $\bigH$ coincides with $\phi_j$, which is a unital completely isometric map. Hence, $\iota_j$ is completely isometric. 

The proof of Theorem \ref{Tm:amalgamated_ucc} lends itself available to obtain a complete characterization of maximal representations on free products. The following forward implication was essentially given by Duncan in \cite[Section 3, Theorem 1]{Duncan_FreeProd}, using multiplicative domain arguments. A gap in the proof of this theorem was filled by Davidson, Fuller and Kakariadis in \cite[Theorem 5.3.20]{Dav-Ful-Kak}, where they show that the free product of operator algebras embeds completely isometrically in the free product of $C^*$-envelopes. We provide a full characterization of maximal representations that yields these results as easy consequences.

\begin{Tm} \label{tm:UEP-preservation-free-prod}
Let $\{\bigA_i\}_{i\in I}$ be a family of unital operator algebras containing a common unital $C^*$-algebra $\bigD$ with $1_{\bigA_i} = 1_{\bigD}$. Suppose further that $\pi_i : \bigA_i \rightarrow B(\bigH)$ are representations that agree on $\bigD$ and denote $\pi := \underset{\bigD}{*} \pi_i$. Then each $\pi_i$ is maximal if and only if the representation $\pi$ is maximal.
\end{Tm}

\begin{proof}
Without loss of generality we may assume all of our representations are unital. Suppose that each $\pi_i$ is maximal, and let $\psi$ be representation dilating $\pi= \underset{\bigD}{*} \pi_i$. Then $\psi|_{\bigA_i}$ dilates each $\pi_i$, and by maximality of $\pi_i$, the subspace $\bigH$ is reducing for each $\psi|_{\bigA_i}$. Since $\psi = \underset{\bigD}{*}(\psi|_{\bigA_i})$, we see that $\bigH$ is reducing for $\psi$. Hence, $\pi$ is maximal.

Conversely, pick $0\in I$, and let $\varphi_0 : \bigA_0 \rightarrow B(\bigK)$ be a representation that dilates $\pi_0$. By Sarason's theorem \cite[Exercise 7.6]{PauBook}, we see that $\bigH$ is \emph{semi-invariant} in the sense that we may write $\bigK = \bigL \oplus \bigH \oplus \bigR$ such that $\bigL$ and $\bigL \oplus \bigH$ are invariant for $\varphi_0(\bigA_0)$. Hence for $a\in \bigA_0$ we have the following block triangular form for $\varphi_0(a)$,
\[
\varphi_0(a) = \begin{bmatrix}
* & * & *\\
0 & \pi_0(a) & * \\
0 & 0 & *
\end{bmatrix}.
\]
Since $\bigD$ is self adjoint, we see that $\bigR,\bigH$ and $\bigL$ are reducing for $\varphi_0(\bigD)$. For $0\neq j \in I$, let $\psi_j : \bigA_j \rightarrow B(\bigL)\oplus B(\bigR)$ be a UCC extension of $d \mapsto \varphi_0(d)|_{\bigL} \oplus \varphi_0(d)|_{\bigR}$.
For $0\neq j \in I$ we set $\phi_j = \pi_j \oplus \psi_j$, and $\phi_0 = \varphi_0$, so that the family $\{\phi_i\}_{i\in I}$ are $\bigD$-bimodule UCC maps that agree on $\bigD$.
We denote $\bigL_1 := \bigL$, $\bigR_1:=\bigR$ and $\phi_i^{(1)}:= \phi_i$ for $i\in I$. From the Arveson-Stinespring's dilation theorem, for each $0\neq j \in I$ we can dilate each $\phi_j$ to a representation $\rho_j^{(1)} : \bigA_j \rightarrow B(\bigL_1^{(j)} \oplus \bigL_1) \oplus B(\bigH) \oplus B(\bigR_1^{(j)} \oplus \bigR_1)$. We also set $\rho_0^{(1)}:= \phi_0 = \varphi_0$. Next, let $n \in \mathbb{N}$, and suppose that 
\begin{enumerate}
\item
For $1 \leq m < n$ and any $i\in I$ we have subspaces $\bigL_m \subseteq \bigL_m \oplus \bigL_m^{(i)} \subseteq \bigL_{m+1}$ and $\bigR_m \subseteq \bigR_m \oplus \bigR_m^{(i)} \subseteq \bigR_{m+1}$.
\item
For every $1\leq m \leq n$ we have unital completely contractive $\bigD$-bimodule maps $\phi_j^{(m)} : \bigA_j \rightarrow B(\bigL_m) \oplus B(\bigH) \oplus B(\bigR_m)$ for $0\neq j \in I$ and $\phi_0^{(m)} : \bigA_0 \rightarrow B(\bigL^{(m)}_0 \ominus \bigL_1) \oplus B(\bigK) \oplus B(\bigR^{(m)}_0 \ominus \bigR_1)$ such that $\{\phi_i^{(m)} \}_{i\in I}$ all agree on $\bigD$.
\item
For every $1\leq m \leq n$ we have representations $\rho_j^{(m)} : \bigA_j \rightarrow B(\bigL_m^{(j)} \oplus \bigL_m) \oplus B(\bigH) \oplus B(\bigR_m^{(j)}\oplus \bigR_m)$ for $0 \neq j \in I$ and $\rho_0^{(m)} : \bigA_0 \rightarrow B(\bigL_m^{(0)}\oplus (\bigL_m \ominus \bigL_1)) \oplus B(\bigK) \oplus B(\bigR_m^{(0)} \oplus(\bigR_m \ominus \bigR_1))$ such that $\rho^{(m)}_i$ dilates $\phi_i^{(m)}$, and both $\rho_i^{(m+1)}$ and $\phi_i^{(m+1)}$ have $\rho_i^{(m)}$ as a direct summand for each $i\in I$.
\item
For every $1\leq m \leq n$ and $i\in I$ the restriction $\phi^{(m)}_i|_{\bigD}$ has $\rho^{(m)}_i|_{\bigD}$ as a direct summand.

\end{enumerate}

Denote $\bigL_{n+1} = \bigL_n \oplus \bigoplus_{i\in I} \bigL_n^{(i)}$ and $\bigR_{n+1} = \bigR_n \oplus \bigoplus_{i \in I} \bigR_n^{(i)}$. For each $i \in I$ and $k\in I$ distinct from $i$, we use Arveson's extension theorem to extend $d \mapsto \rho_k^{(n)}(d)|_{\bigL_n^{(k)}} \oplus \rho_k^{(n)}(d)|_{\bigR_n^{(k)}}$ to a UCC $\bigD$-bimodule map $\sigma_{ik} : \bigA_i \rightarrow B(\bigL_n^{(k)}) \oplus B(\bigR_n^{(k)})$. We then define for all $i \in I$ the UCC $\bigD$-bimodule maps
\begin{equation} \label{eq:phi-j}
\phi_i^{(n+1)} = \rho_i^{(n)} \oplus \bigoplus_{k\neq i}\sigma_{ik}
\end{equation}
so that for $0\neq j \in I$ we have $\phi_j^{(n+1)} : \bigA_j \rightarrow B(\bigL_{n+1}) \oplus B(\bigH) \oplus B(\bigR_{n+1})$, and $\phi_0^{(n+1)}: \bigA_0 \rightarrow B(\bigL_{n+1} \ominus \bigL_1) \oplus B(\bigK) \oplus B(\bigR_{n+1} \ominus \bigR_1)$.
Hence, we see for each $i\in I$ that $\rho_i^{(n)}$ is a direct summand of $\phi_i^{(n+1)}$ and for $d\in \bigD$ that
$$
\phi_i^{(n+1)}(d) = \phi_i^{(n)}(d) \oplus \bigoplus_{k\in I} (\rho_k^{(n)}(d)|_{\bigL_n^{(k)}} \oplus \rho_k^{(n)}(d)|_{\bigR_n^{(k)}}).
$$
Thus, $\{\phi_i^{(n+1)}\}_{i\in I}$ agree on $\bigD$, since all $\{\phi_i^{(n)}\}_{i\in I}$ agree on $\bigD$.

We may now use the special form of $\phi_i^{(n+1)}$ in equation \eqref{eq:phi-j} to dilate each $\phi_j^{(n+1)}$ for $0\neq j \in I$ to a representation $\rho_j^{(n+1)} : \bigA_j \rightarrow B(\bigL_{n+1}^{(j)} \oplus \bigL_{n+1}) \oplus B(\bigH) \oplus B(\bigR_{n+1}^{(j)} \oplus \bigR_{n+1})$, along with $\rho_0^{(n+1)} : \bigA_0 \rightarrow B(\bigL_{n+1}^{(0)} \oplus (\bigL_{n+1} \ominus \bigL_1)) \oplus B(\bigK) \oplus B(\bigR_{n+1}^{(0)} \oplus (\bigR_{n+1} \ominus \bigR_1))$ dilating $\phi_0^{(n+1)}$. So that now $\rho_i^{(n)}$ is a direct summand of $\rho_i^{(n+1)}$ for each $i\in I$.

After constructing an infinite family $\phi_i^{(m)}$ and $\rho_i^{(m)}$ for any $m\in \mathbb{N}$ satisfying properties $(1) - (4)$ we proceed as follows. We define for $0 \neq j \in I$ representations $\rho_j : \bigA_j \rightarrow B(\bigL') \oplus B(\bigH) \oplus B(\bigR')$ where $\bigL' = \bigvee_{n \geq 1} \bigL_n$ and $\bigR' = \bigvee_{n \geq 1} \bigR_n$ by setting $\rho_j(a)h = \rho_j^{(n)}(a)h$ for $h \in \bigL_n \oplus \bigL_n^{(j)} \oplus \bigH \oplus \bigR_n \oplus \bigR_n^{(j)}$. And similarly define a representation $\rho_0 : \bigA_0 \rightarrow B(\bigL' \ominus \bigL_1) \oplus B(\bigK) \oplus B(\bigR' \ominus \bigR_1)$.

Then as in the proof of Theorem \ref{Tm:amalgamated_ucc} we see that $\{\rho_i \}_{i\in I}$ are representation on the Hilbert space $\bigL' \oplus \bigH \oplus \bigR'$ that agree on $\bigD$.

Now, since $\bigH$, $\bigL'$ and $\bigR'$ are reducing for all $\rho_j$ with $0 \neq j \in I$, and since $\bigH$ is semi-invariant for $\rho_0$, while $\bigL' \ominus \bigL_1$ and $\bigR' \ominus \bigR_1$ are reducing for $\rho_0$, we see that the compression of $\rho:= \underset{\bigD}{*}\rho_i$ to $\bigH$ is $\pi$. Indeed, by induction on $\ell$ we can show for $a_{i_1} \in \bigA_{i_1}, ..., a_{i_n}\in \bigA_{i_{\ell}}$ elements with $i_{m} \neq i_{m+1}$ for $1 \leq m <\ell$, then
$$
\rho(a_{i_1} * ... * a_{i_{\ell}})P_{\bigH} = \prod_{m=1}^{\ell} P_{\bigL' \oplus\bigH}\rho_{i_m}(a_{i_m})P_{\bigL' \oplus\bigH}
$$
and
$$
P_{\bigH} \rho(a_{i_1} * ... * a_{i_{\ell}}) = \prod_{m=1}^{\ell} P_{\bigH \oplus \bigR'}\rho_{i_m}(a_{i_m})P_{\bigH \oplus \bigR'}
$$
So that
$$
P_{\bigH} \rho(a_{i_1} * ... * a_{i_{\ell}})P_{\bigH} = \prod_{m=1}^{\ell} P_{\bigH}\rho_{i_m}(a_{i_m})P_{\bigH} =
$$
$$
\prod_{m=1}^{\ell} \pi_{i_m}(a_{i_m}) = \pi(a_{i_1} * ... * a_{i_{\ell}}).
$$
Thus, we see that $\rho$ is a dilation of $\pi$. But since $\pi$ is maximal, we see that $\bigH$ is reducing for $\rho$, and in particular reduces $\rho_0$. However, $\bigK$ reduces $\rho_0$ and $\bigH \subseteq \bigK$, so that $\bigH$ reduces $\rho_0 = \phi_0 = \varphi_0$. Hence, we see that $\pi_0$ is maximal as well. Since $0\in I$ was arbitrary, we see that all $\pi_i$ are maximal.
\end{proof}

In \cite[Section 4]{BlecherPaulsen_Explicit}, Blecher and Paulsen prove the complete injectivity of the free product of operator algebras amalgamated over the complex numbers. We next prove this where the amalgamation is over any common $C^*$-algebra. This generalizes \cite[Proposition 2.2]{ADER-free-prod} due to Armstrong, Dykema, Exel and Li for free products of finitely many $C^*$-algebras amalgamated over a common $C^*$-algebra. 

We note here that the complete isometric embedding of the free product of operator algebras inside the free product of their $C^*$-envelopes was first proven by Davidson, Fuller and Kakariadis in \cite[Theorem 5.3.20]{Dav-Ful-Kak}, and that their methods can be adapted to prove complete injectivity in our context. We will instead provide a proof of this fact by using maximality.

\begin{Pn}\label{Pn:CompleteInjectivity}
The free product of unital operator algebras amalgamated over a common $C^*$-subalgebra is completely injective. That is, if $\{\bigA_i\}_{i\in I}$ and $\{\bigB_i\}_{i\in I}$ are two families of unital operator algebras containing a common $C^*$-subalgebra $\bigD$ such that $\bigA_i$ is an operator subalgebra of $\bigB_i$ for every $i \in I$ and $1_{\bigA_i} = 1_{\bigB_i} = 1_{\bigD}$, then the inclusion $\underset{\bigD}{\hat{*}} \bigA_i\subseteq \underset{\bigD}{\hat{*}} \bigB_i$ is completely isometric. 
\end{Pn}

\begin{proof}
Denote $\hat \bigA :=\underset{\bigD}{\hat{*}} \bigA_i$ and $\hat \bigB : = \underset{\bigD}{\hat{*}} \bigB_i$.
Let $\iota_i:\bigB_i \to \hat\bigB$ and $\kappa_i:\bigA_i \to \hat\bigA$ denote the natural completely isometric inclusions. Then the unital completely isometric homomorphisms $\iota_i|_{\bigA_i}:\bigA_i \to  \hat\bigB$ agree on $\bigD$, so $\phi:=\underset{\bigD}{*} (\iota_i|_{\bigA_i}):\hat\bigA \to \hat\bigB$ is a unital completely contractive homomorphism. If we denote by $\| \cdot \|_{\hat{\bigB}, n}$ and $\| \cdot \|_{\hat{\bigA},n}$ the matricial $n$-norms on $\hat{\bigA}$ and $\hat{\bigB}$ respectively, it is clear that $\|A\|_{\hat \bigA, n} \leq \|A\|_{\hat \bigB, n}$ for any $A\in M_n(\hat \bigA)$.

Hence, it would suffice to show that $\| A \|_{\hat{\bigA},n} \leq \| A \|_{\hat{\bigB},n}$ for every $A \in M_n(\hat{\bigA})$. To this end, identify a completely isometric and \emph{maximal} copy of $\hat \bigA$ inside $B(\bigH)$ for some Hilbert space $\bigH$. By Theorem \ref{tm:UEP-preservation-free-prod} we see that $\kappa_i:\bigA_i \to\hat \bigA \subseteq B(\bigH)$ are maximal. By Arveson's extension theorem the maps $\kappa_i:\bigA_i \to\hat \bigA \subseteq B(\bigH)$ extend to unital completely contractive maps $\tilde \kappa_i:\bigB_i \to B(\bigH)$ which agree on $\bigD$, and by (the proof of) Theorem \ref{Tm:amalgamated_ucc}, there are representations $\rho_i : \bigB_i \rightarrow B(\bigK)$ that dilate each $\widetilde{\kappa_i}$ and agree on $\bigD$. Then, each restriction $\rho_i|_{\bigA_i}$ is a dilation of $\kappa_i$, and as each $\kappa_i$ is maximal, $\bigH$ must reduce $\rho_i|_{\bigA_i}$. We denote the free product by $\rho = \underset{\bigD}{*}\rho_i$. Thus, for each $A\in M_n(\hat{\bigA})$ we have that 
\[
P_{\bigH}\rho(A)|_{\bigH} = \underset{\bigD}{*}P_{\bigH}\rho_i(A)|_{\bigH} = \underset{\bigD}{*} \kappa_i(A) = A.
\]
Hence, for $A \in M_n(\hat{\bigA})$, we have
\[
\|A\|_{\hat{\bigA},n} = \|P_{\bigH}\rho(A)|_{\bigH}\| \leq \|\rho(A)\| \leq \|A\|_{\hat{\bigB},n},
\]
and the proof is complete.
\end{proof}

When $\{\bigA_i\}_{i\in I}$ are all non-unital with a common non-unital $C^*$-algebra $\bigD$, we define their free product $\underset{\bigD}{\hat{*}} \bigA_i$ to be the operator algebra generated by the images of $\bigA_i$ inside the free product of their unitization $\underset{\bigD^1}{\hat *} \bigA_i^1$ amalgamated over the unitization $\bigD^1$. By Meyer's theorem and the proof of \cite[Lemma 2.3]{ADER-free-prod}, we similarly get that $\underset{\bigD^1}{\hat *} \bigA_i^1$ coincides with the unitization $(\underset{\bigD}{\hat{*}} \bigA_i)^1$. Using this, it follows that the non-unital free product shares the analogous pushout universal property and complete injectivity. By the discussion together with Proposition \ref{P:unital-red}, we obtain a not-necessarily-unital version of Theorem \ref{tm:UEP-preservation-free-prod}.

\begin{Cy} \label{cor:non-unital-pres}
Let $\{\bigA_i\}_{i\in I}$ be a family of \emph{either all unital, or all non-unital} operator algebras, each containing a common $C^*$-algebra $\bigD$. In the unital case, assume further that $\bigD$ is unital with a common unit. Suppose further that $\pi_i : \bigA_i \rightarrow B(\bigH)$ are representations that agree on $\bigD$, and denote $\pi := \underset{\bigD}{*} \pi_i$. Then each $\pi_i$ is maximal if and only if $\pi$ is maximal.
\end{Cy}

Using complete injectivity, given a family of operator algebras $\{\bigA_i\}_{i \in I}$ with a common $C^*$-subalgebra $\bigD$, we will henceforth freely identify $\underset{\bigD}{\hat{*}} \bigA_i$ as a subalgebra of $\underset{\bigD}{\hat{*}} \bigB_i$, where $\bigB_i$ is any operator algebra containing $\bigA_i$. Now that the basic theory of free products has been established, we will abuse notation and denote $\underset{\bigD}{*} \bigA_i$ for the operator algebraic free product $\underset{\bigD}{\hat{*}} \bigA_i$.

\begin{Rk} \em
In a previous version of this manuscript, Corollary \ref{cor:non-unital-pres} (which used to be Proposition 4.4 in the previous version) was proven using Boca's theorem \cite[Theorem 3.1]{Boca}, and under the additional assumption that each $\bigA_i$ generate $C^*$-algebras $\bigB_i$ that have conditional expectations $E_i : \bigB_i \rightarrow \bigD$, such that $\bigA_i = \bigD \oplus (\bigA_i \cap \ker E_i)$ as a $\bigD$-bimodule. It is important to mention here that in \cite{Dav-Kak}, a short proof of Boca's theorem is given using dilation techniques similar to ours. This justifies the point of view that the methods used here and in \cite{Dav-Kak, Dav-Ful-Kak} are dilation theoretic \emph{generalizations} of Boca's theorem, that no longer require conditional expectations. We are grateful to Elias Katsoulis for pointing out an issue in the previous proof of Proposition \ref{Pn:CompleteInjectivity}, which led to the full characterization of (what is now) Corollary \ref{cor:non-unital-pres}.
\end{Rk}

\section{Full Cuntz-Krieger dilations for free families} \label{sec:free-prod-graph}

In this section we will write $G=(V,E)$ for a directed graph, where the source and range maps are understood implicitly. A function $c: E \rightarrow I$ is an $I$-coloring of edges of $G$, and we define an $I$-colored graph to be the triple $(V,E,c)$. Given such a colored directed graph $(V,E,c)$, we denote $E_i = c^{-1}(i)$.

Let $G=(V,E,c)$ be an $I$-colored directed graph and $\bigH$ a Hilbert space. A {\em Toeplitz-Cuntz-Krieger $I$-colored family} on $\bigH$ is a pair $(P,S)$ comprised of a operators $P:= \{P_v : v \in V\}$ on $H$ and an $I$-tuple of sets of operators $S:= \{S^{(i)}\}_{i=1}^d$ on $\bigH$ with $S^{(i)}:= \{S^{(i)}_e : e \in c^{-1}(i) \}$ such that each $(P,S^{(i)})$ is a a TCK family for $(V,E_i)$ for each $i\in I$. 

We say that $(P,S)$ is a {\em Cuntz-Krieger $I$-colored family} / {\em full Cuntz-Krieger $I$-colored family} if, in addition, each $(P,S^{(i)})$ is a a CK / full CK family for $(V,E_i)$ for each $i\in I$ respectively.

From here on out, for a given set of vertices $V$, we set $\bigV:=C_0(V)$. We will identify $\bigV$ with $C^*(\{P_v\})$ for some (all) TCK or CK $I$-colored families $(P,S)$ where $P_v\neq 0$ for all $v\in V$ and any colored graph $G$ with $V$ as its set of vertices.

Given a colored directed graph $G=(V,E,c)$, let $G_i = (V,c^{-1}(i))$ be the graph of color $i \in I$. By compounding universal properties, it is easy to see that TCK $I$-colored families are in bijection with $*$-representations of the free product $\underset{\bigV}{*} \bigT(G_i)$ over $i\in I$ and that CK $I$-colored families are in bijection with $*$-representations of the free product $\underset{\bigV}{*} \bigO(G_i)$ over $i\in I$. We will call a $*$-representation of either $\underset{\bigV}{*}\bigO(G_i)$ of $\underset{\bigV}{*}\bigT(G_i)$ full Cuntz-Krieger if its associated TCK $I$-colored family is full CK $I$-colored family.

We next characterize those representations of the free product of Toeplitz graph algebras with the unique extension property when restricted to the free product of graph tensor algebras. Recall that for countable directed graphs $\{G_i\}_{i\in I}$ on the same vertex set $V$, by Proposition \ref{Pn:CompleteInjectivity}, we have that the embedding of $\underset{\bigV}{*} \bigT_+(G_i)$ in $\underset{\bigV}{*} \bigT(G_i)$ is completely isometric, so we may identify the former as a subalgebra of the latter.

\begin{Pn} \label{prop:free-prod-uep-char}
Let $\{G_i\}_{i\in I}$ be a collection of countable directed graphs on the same vertex set $V$, and let $\pi : \underset{\bigV}{*} \bigT(G_i) \rightarrow B(\bigH)$ be a $*$-representation. Then $\pi|_{\underset{\bigV}{*} \bigT_+(G_i)}$ has the unique extension property if and only if for each $i\in I$ the $*$-representation $\pi_i : = \pi |_{\bigT(G_i)}$ is full CK with respect to $G_i$.
\end{Pn}

\begin{proof}
By Corollary \ref{cor:non-unital-pres} we see that $\pi : \underset{\bigV}{*} \bigT(G_i) \rightarrow B(\bigH)$ has the unique extension property when restricted to $\underset{\bigV}{*}\bigT_+(G_i)$ if and only if each $\pi_i$ has the unique extension property when restricted to $\bigT_+(G_i)$. By Theorem \ref{Tm:non-annihi-no-uep}, this occurs if and only if each $\pi_i$ is full CK with respect to $G_i$.
\end{proof}

We apply Proposition \ref{prop:free-prod-uep-char} to draw a dilation result that generalizes Corollary \ref{Cy:full-CK-dilation}.

\begin{Cy} \label{Cy:colored-dilation}
Let $G=(V,E,c)$ be an $I$-colored directed graph, and let $(P,S)$ be a TCK $I$-colored family on $\bigH$. Then there exists a full CK $I$-colored family $(Q,T)$ on a Hilbert space $K$ containing $\bigH$, such that $f(P,S) = P_{\bigH} f(Q,T)|_{\bigH}$ for any polynomial $f \in \mathbb{C}\langle V,E \rangle$ in non-commuting variables.
\end{Cy}

\begin{proof}
Let $\pi_{P,S} : \underset{\bigV}{*} \bigT(G_i) \rightarrow B(\bigH)$ be the $*$-representation associated to $(P,S)$. By \cite[Theorem 1.2]{DritschelMcCullough} we can dilate $\pi_{P,S}|_{\underset{\bigV}{*} \bigT_+(G_i)}$ to a maximal completely contractive homomorphism $\tau : \underset{\bigV}{*} \bigT_+(G_i) \rightarrow B(\bigK)$. Without loss of generality, $\bigH$ is a subspace of $\bigK$. Let $\rho : \underset{\bigV}{*} \bigT(G_i) \rightarrow B(\bigK)$ be its unique extension to a $*$-representation, and $(Q,T)$ the associated TCK family of $\rho$. As $\tau$ has the unique extension property, by Proposition \ref{prop:free-prod-uep-char} we must have that $(Q,T)$ is full CK, and as $\tau$ dilates $\pi_{P,S}|_{\underset{\bigV}{*} \bigT_+(G_i)}$, we have that every polynomial $f \in \mathbb{C}\langle V,E \rangle$ in non-commuting variables must satisfy $f(P,S) = P_{\bigH} f(Q,T)|_{\bigH}$.
\end{proof}

The following result mirrors \cite[Proposition 1.6]{Skalski-Zacharias-Wold} on the existence of maximal fully co-isometric summands, but when restricting to the context of directed graphs our result is more general as it requires no relations between families of different colors.

\begin{Cy}
Let $G = (V,E,c)$ be an $I$-colored directed graph, and let $(P,S)$ be a TCK $I$-colored family on $\bigH$. Then there is a unique maximal common reducing subspace $\bigK$ for operators in $(P,S)$ such that $(P,S)|_{\bigK}$ is full CK.
\end{Cy}

\begin{proof}
Let $\pi_{P,S}$ be the $*$-representation associated to $(P,S)$. By Proposition \ref{prop:largest-UEP} there is a unique largest reducing subspace $\bigK$ for $\pi_{P,S}$ such that $\rho: \underset{\bigV}{*} \bigT(G_i) \rightarrow B(\bigK)$ given by $\rho(b) = \pi_{P,S}(b)|_{\bigK}$ has the unique extension property when restricted to $\underset{\bigV}{*} \bigT_+(G_i)$. By Proposition \ref{prop:free-prod-uep-char}, we see that the associated TCK family $(P,S)|_{\bigK}$ is in fact full CK, and $\bigK$ is a unique maximal common reducing subspace with this property.
\end{proof}

Denote by $\bigH_V = \bigoplus_{v\in V} \bigH_v$ a Hilbert space direct sum of separable infinite dimensional Hilbert spaces $\bigH_v$. Then $\bigV \cong C^*(\{p_v\})$ can be represented on $\bigH_V$ by mapping $p_v$ to the projection $P_v$ onto $\bigH_v$, and if $\rho : \bigT(G_i) \rightarrow B(\bigH)$ is a non-degenerate representation such that $\rho|_{C^*(\{p_v\})} : \bigV \cong C^*(\{p_v\}) \rightarrow B(\bigH)$ where $\rho(p_v)\bigH$ is infinite dimensional for each $v\in V$, then $\rho$ is unitarily equivalent to a representation on $\bigH_V$ where $\rho(p_v) = P_v$ is the projection onto $\bigH_v$.

When all $G_i$ are row-finite, Duncan \cite[Proposition 4.4]{Duncan_EdgeColored} explained how $\underset{\bigV}{*} \bigT_+(G_i)$ has the unique extension property in $\underset{\bigV}{*} \bigO(G_i)$. This result also relied on \cite[Section 3, Theorem 1]{Duncan_FreeProd}, and as mentioned earlier, a gap in that result was later filled in \cite[Theorem 5.3.20]{Dav-Ful-Kak}. We next prove Duncan's graph-algebra result while providing its converse.

\begin{Tm}\label{Tm:freeHR}
Let $\{G_i\}_{i\in I}$ be a collection of countable directed graphs over the same vertex set $V$. Then the quotient map $q:\underset{\bigV}{*} \bigT(G_i) \rightarrow \underset{\bigV}{*} \bigO(G_i)$ is completely isometric on $\underset{\bigV}{*} \bigT_+(G_i)$. Hence $C_e^*(\underset{\bigV}{*} \bigT_+(G_i))$ is a quotient of $\underset{\bigV}{*} \bigO(G_i)$.

Moreover, we have that each $G_i$ is row-finite if and only if $\underset{\bigV}{*} \bigT_+(G_i)$ has the unique extension property in $\underset{\bigV}{*} \bigO(G_i)$. In particular, in this case we have $C_e^*(\underset{\bigV}{*} \bigT_+(G_i)) \cong \underset{\bigV}{*} \bigO(G_i)$.
\end{Tm}

\begin{proof}
Since $\bigT_+(G_i)$ can be identified as a subalgebra of $\bigO(G_i)$ via the image of the quotient map $q_i: \bigT(G_i) \rightarrow \bigO(G_i)$, by Proposition \ref{Pn:CompleteInjectivity}, we see that $\underset{\bigV}{*} \bigT_+(G_i)$ is can be identified as a subalgebra of $\underset{\bigV}{*} \bigO(G_i)$ via the image of $q = \underset{\bigV}{*}q_i$.

For the second part, suppose each $G_i$ is row-finite. Let $\pi : \underset{\bigV}{*}\bigO(G_i) \rightarrow B(\bigH)$ be a $*$-representation. Then $\pi \circ q$ is a $*$-representation of $\underset{\bigV}{*}\bigT(G_i)$, and by invariance of the UEP it will suffice to show that $(\pi \circ q)|_{\underset{\bigV}{*}\bigT_+(G_i)}$ has the UEP. By Proposition \ref{prop:free-prod-uep-char}, this happens if and only if $\pi_i \circ q_i$ is full CK. However, as each $G_i$ is row-finite, Theorem \ref{Tm:non-annihi-no-uep} implies that each $\pi_i \circ q_i : \bigT(G_i) \rightarrow B(\bigH)$ is full CK. Hence, $\underset{\bigV}{*} \bigT_+(G_i)$ has the unique extension property in $\underset{\bigV}{*} \bigO(G_i)$.

For the converse, we will construct a $*$-representation of $\underset{\bigV}{*} \bigO(G_i)$ that lacks the unique extension property when restricted to $\underset{\bigV}{*} \bigT_+(G_i)$. Indeed, as one of $\{G_i\}$ is not row-finite, by Theorem \ref{theorem:C-envelope-hyperrigidity} there is some $j\in I$ for which there is a CK representation $\rho_j : \bigT(G_j) \rightarrow B(\bigH)$ that is not a full CK representation, and up to inflating $\bigH$ we may assume $\rho_j = \rho_j^{(\infty)}$. For $i\in I$ different from $j$, let $\rho_i : \bigT(G_i) \rightarrow B(\bigH)$ be any representation annihilating the Cuntz-Krieger ideal $\bigJ(G_i)$ for which $\rho_i(p_v) \neq 0$ for all $v\in V$. Again up to inflating $\bigH$ we may assume $\rho_i = \rho_i^{(\infty)}$. In this case for all $i\in I$ the representation $\rho_i$ is unitarily equivalent to a representation on $\bigH_V$ where $\rho_i(p_v)$ is mapped to the projection $P_v$. In this case the free product $\underset{\bigV}{*}\rho_i$ is well-defined, and by Proposition \ref{prop:free-prod-uep-char} we see that $\underset{\bigV}{*}\rho_i$ does not have the unique extension property when restricted to $\underset{\bigV}{*}\bigT_+(G_i)$ while still annihilating $\bigJ = \langle \bigJ(G_i)\rangle_{i \in I}$, so that it induces a representation $\underset{\bigV}{*}\dot{\rho_i}$ on $\underset{\bigV}{*}\bigO(G_i)$ that does not have the unique extension property.
\end{proof}

\begin{Rk}\em
After the completion of this work we were informed by Evgenios Kakariadis, that it is possible to show that $C_e^*(\underset{\bigV}{*} \bigT_+(G_i)) \cong \underset{\bigV}{*} \bigO(G_i)$, regardless of whether the graphs $G_i$ are row-finite or not. This completes the picture described in Theorem \ref{Tm:freeHR}.
\end{Rk}

\section*{Acknowledgments}
The first author wishes to thank Vern Paulsen and Samuel Harris for illuminating discussions on free products. Both authors wish to thank their PhD advisors, Ken Davidson and Orr Shalit, for reading and providing remarks on a draft version of this paper. A special thanks goes to the anonymous referee, whose remarks led to an improved version of Section \ref{sec:ncb}.


\end{document}